\documentclass{elsart}
\journal{Journal of Symbolic Computation}
\usepackage{amsmath}
\usepackage{natbib}
\usepackage{yjsco}
\usepackage{amsfonts}
\usepackage{amssymb}
\usepackage{amscd}

\DeclareMathOperator{\rank}{rk}
\DeclareMathOperator{\Char}{char}
\DeclareMathOperator{\Rem}{\sf d-rem}
\DeclareMathOperator{\drem}{\sf d-rem}
\DeclareMathOperator{\ord}{ord}
\DeclareMathOperator{\charact}{char}

\newcommand{\G}{Gr\"obner}

\newcommand{\K}{{\mathbf k}}

\newcommand{\C}{{\mathbb C}}
\newcommand{\A}{{\mathbb A}}

\def\B{{\mathbb B}}
\def\p{{\mathfrak p}}
\def\q{{\mathfrak q}}
\def\D{{\mathbb D}}

\newcommand{\ld}{\mathop{\rm ld}\nolimits}

\newcommand{\lv}{\mathop{\rm lv}\nolimits}

\newcommand{\rk}{\mathop{\rm rk}\nolimits}

\newcommand{\algrem}{\mathop{\sf algrem}\nolimits}

\newcommand{\sep}{{\;|\;}}
                                                 
\renewcommand{\u}{{\bf u}}
\renewcommand{\i}{{\bf i}}
\newcommand{\s}{{\bf s}}

\renewcommand{\b}{\ \ \ }

\newcommand{\Le}{\leqslant}
\newcommand{\Ge}{\geqslant}

\begin{document}

\begin{frontmatter}
\title{Bounds for algorithms in differential algebra\thanksref{thank1}}
\author[OG]{Oleg Golubitsky\thanksref{thank3}}
\ead{oleg.golubitsky@gmail.com}
\ead[url]{http://publish.uwo.ca/\textasciitilde ogolubit/}
\author[MVK]{Marina Kondratieva}
\ead{kondrmar@rol.ru}
\ead[url]{http://shade.msu.ru/\textasciitilde kondra\_m/}
\author[MMM]{Marc Moreno Maza\thanksref{thank4}}
\ead{moreno@orcca.on.ca}
\ead[url]{http://www.csd.uwo.ca/\textasciitilde moreno/}
\author[AO]{Alexey Ovchinnikov\thanksref{thank2}}
\address[OG]{University of Western Ontario\\Department of Computer Science\\London, Ontario, Canada N6A 5B7}
\address[MVK]{Moscow State University\\ Department of Mechanics and Mathematics\\ Leninskie gory, Moscow, Russia, 119992}
\address[MMM]{University of Western Ontario\\Department of Computer Science\\London, Ontario, Canada N6A 5B7}
\address[AO]{North Carolina State University\\ Department of Mathematics\\ Raleigh, NC 27695-8205, USA}
\date\today
\ead{aiovchin@ncsu.edu}
\ead[url]{http://www4.ncsu.edu/\textasciitilde aiovchin/}
\thanks[thank1]{The work was partially supported by the Russian Foundation for Basic Research, project no. 05-01-00671.}
\thanks[thank2]{This author was also partially supported by NSF Grant CCR-0096842.}
\thanks[thank3]{This author was also partially supported by NSERC Grant PDF-301108-2004.}
\thanks[thank4]{This author was also partially supported by NSERC Grant RGPIN     {\em Algorithms and software for triangular decompositions of algebraic and differential systems.}}

\begin{keyword}
differential algebra \sep characteristic sets \sep radical differential ideals \sep decomposition into regular components
\MSC 12H05 \sep 13N10 \sep 13P10
\end{keyword}

\begin{abstract} 
  We consider the Rosenfeld-Gr\"obner algorithm for computing a regular 
  decomposition of a radical differential ideal generated by a set
  of ordinary differential polynomials in $n$ indeterminates. 
  For a set of ordinary differential polynomials $F$, let $M(F)$ be the sum 
  of maximal orders of differential indeterminates occurring in $F$.  
  We propose a modification of the Rosenfeld-Gr\"obner algorithm, in which 
  for every intermediate polynomial system $F$, the bound
  $M(F)\Le (n-1)!M(F_0)$ holds, where $F_0$ is the initial set of generators
  of the radical ideal. In particular, the resulting regular systems satisfy the bound. 
  Since regular ideals can be decomposed into characterizable
  components algebraically, the bound also holds
  for the orders of derivatives occurring in a characteristic
  decomposition of a radical differential ideal.

  We also give an algorithm for converting a characteristic decomposition
  of a radical differential ideal from one ranking into another. This
  algorithm performs all differentiations in the beginning
  and then uses a purely algebraic decomposition algorithm.
\end{abstract}

\end{frontmatter}

\section{Introduction}
This paper is about constructive differential algebra. We study algorithms
dealing with algebraic differential equations. Many different problems
can be addressed to this topic. One can, for instance, test membership
to a radical differential ideal, compute the Kolchin dimensional polynomial.
The kind of algorithms we are dealing with are decomposition algorithms for 
radical differential ideals. Generally, there are two such algorithms, although
they have variations. 

The Ritt-Kolchin algorithm  
computes a prime decomposition of a radical differential ideal, where each prime
component is represented by its characteristic set. This algorithm is based
on important results in differential algebra (see \cite{Rit,Kol}), such as the Basis Theorem,
the Prime Decomposition Theorem for radical differential ideals, 
the differential version of the Hilbert Theorem of Zeroes, 
and the Rosenfeld Lemma. 
It also relies on the solution of the so-called factorization problem: 
given an autoreduced set, determine whether the corresponding algebraic saturated ideal 
is prime and, if it is not, find two polynomials outside of the ideal whose product belongs 
to the ideal. 

Due to the complexity of the factorization problem, it was desirable to 
avoid it, which was done in the Rosenfeld-Gr\"obner algorithm proposed in
\citep{Bou1}. Instead of decomposing a given radical differential ideal
into prime components, this algorithm represents it as an intersection of regular
differential ideals, also introduced in \citep{Bou1}; the correctness of the algorithm, 
in addition to the above-mentioned
theorems, is provided by the Lazard Lemma, which states that regular ideals are radical.
Different proofs of this lemma can be found in \citep{Bou2,Mor,Fac,BLM06}. 

The Rosenfeld-Gr\"obner algorithm is the first decomposition algorithm in differential
algebra that has been actually implemented upto our knowledge. It forms an integral part of the {\tt diffalg} 
package in the computer algebra system Maple. Updates of this package are available 
at {\tt http://www-sop.inria.fr/cafe/Evelyne.Hubert/diffalg/}.
A more efficient implementation of this algorithm in C language 
can be found at the website {\tt http://www.lifl.fr/$\sim$boulier/BLAD/}.

Various improvements of the Rosenfeld-Gr\"obner algorithm have been proposed 
in \citep{Bou2,Fac,Dif,Imp,Und}. They all avoid the factorization problem
and for this reason are called factorization-free methods
in differential algebra. However, no theoretical bound for the computational
complexity of any of these algorithms is known. 

We make the first step towards the goal of estimating this
complexity: we bound the orders of differential polynomials appearing
in the computations. The main results of this work are proven 
only for the {\bf ordinary} case. We consider
the following {\bf two} bounding problems.
The {\bf first} problem is to bound the orders of all intermediate polynomials
and the output of the Rosenfeld-Gr\"obner algorithm. In order to obtain
such a bound in Proposition~\ref{p:finalalgorithm}, we have modified
this algorithm (see Algorithms~\ref{RGBound} and \ref{RGBoundRI}) a little bit.

It would be good to have a bound that would tell us how many
times we need to differentiate the original system in the beginning
of the algorithm, so that the rest of the computation can be performed
by a purely algebraic decomposition algorithm. Since for algebraic
decomposition algorithms complexity estimates are known (see
\cite{Agnes}), such a bound would yield a complexity estimate for 
the differential decomposition as well. In this paper, however,
we do not provide such a bound and, moreover, conjecture that
it would have solved the Ritt problem \citep{Rit}. We leave the 
discovery of such bound and/or the proof of this conjecture
for future research.
%

Nevertheless, for the {\bf second} type of the algorithms we are looking at in this
paper we obtain such a bound. Namely, we can tell how many times one
needs to differentiate elements of a {\it given characteristic set} of a characterizable
differential ideal w.r.t. {\it one} differential ranking, in order
to obtain a characteristic {\it decomposition} of this ideal w.r.t. {\it another}
ranking. In other words, we give a bound for the conversion 
algorithm (Algorithm~\ref{a:conversionalgorithm}) for a characterizable ideal from one ranking
to another (see \cite{Kaehler,PARDI,DGW} for other conversion algorithms applicable
to prime differential ideals).
We emphasize that the input ideal does not have to be characterizable
w.r.t. the target ranking. We show how to obtain its new characteristic
decomposition by first differentiating the input characteristic
set and then applying only algebraic operations (i.e., 
a purely algebraic decomposition algorithm).

The paper is organized as follows. We give an introduction into
differential algebra in Section~\ref{sec:diffalgintro}. Then we describe
the original Rosenfeld-Gr\"obner algorithm in Section~\ref{sec:rg}.
Section~\ref{sec:modifiedRG} is devoted to the bound on the orders
of derivatives computed by a modified version of 
the Rosenfeld-Gr\"obner algorithm. After that, we show how to
transform a characteristic set of a characterizable differential
ideal into a characteristic decomposition of this ideal w.r.t. another
differential ranking. We first do this for prime differential ideals
(Section~\ref{sec:prime}) and then treat the characterizable case in
Section~\ref{sec:characterizable}.

\section{Definitions and notation}\label{sec:diffalgintro}

Differential algebra studies systems of polynomial partial 
differential equations from
the algebraic point of view. The approach is based on the concept 
of differential ring introduced by Ritt.
Recent tutorials on the constructive theory of differential ideals 
are presented in \cite{BouChaptire, BouThese, Dif, Sit}. 
A differential ring is a commutative ring with the unity endowed
with a set of derivations $\Delta = \{\delta_1,\ldots,\delta_m\}$,
which commute pairwise. 
The case of $\Delta = \{\delta\}$ is called {\it ordinary}. If $R$ 
is an ordinary differential ring and $y \in R$, we denote $\delta^ky$
 by $y^{(k)}$.

 Construct the multiplicative monoid 
$\Theta = \left\{\delta_1^{k_1}\delta_2^{k_2}\cdots\delta_m^{k_m}\;\big|\; k_i \Ge 0\right\}$
of {\it derivative operators}. Let $Y=\{y_1,\ldots,y_n\}$ be a set whose
elements are called {\it differential indeterminates}. 
The elements of the set $\Theta Y=\{\theta y\;|\;\theta\in\Theta,\;y\in Y\}$ 
are called {\it derivatives}. Derivative operators from $\Theta$ act on 
derivatives as 
$\theta_1(\theta_2y_i) = (\theta_1\theta_2)y_i$ for all 
$\theta_1,\theta_2 \in \Theta$ and $1 \Le i \Le n$.

The ring of {\it differential polynomials} in differential 
indeterminates $Y$ over a differential field $\K$ is a ring of commutative 
polynomials with coefficients in $\K$ in the infinite set of variables
$\Theta Y$ (see \cite{Kol, Pan, Rit}). 
This ring is denoted $\K\{y_1,\dots,y_n\}$ or $\K\{Y\}$. We consider the case 
of $\charact \K=0$ only. An ideal $I$ in $\K\{Y\}$ is called {\it
  differential}, if for all $f\in I$ and $\delta\in\Delta$, $\delta
f\in I$. We denote differential polynomials by $f, g, h, \ldots$ 
and use letters $I, J, \p$ for ideals.

Let $F \subset \K\{y_1,\dots,y_n\}$ be a set of differential 
polynomials. For the differential and radical differential 
ideal generated by $F$ in $k\{y_1,\dots,y_n\}$, we use 
notations $[F]$ and $\{F\}$, respectively. 

We need the notion of reduction for algorithmic computations.
 First, we introduce a {\it ranking} on the set of derivatives.
A ranking \citep{Kol} is a total order $>$ on the set $\Theta Y$
satisfying the following conditions for all $\theta\in\Theta$ and
$u,v\in\Theta Y$:
\begin{enumerate}
\item $\theta u \Ge u,$
\item $u \Ge v \Longrightarrow \theta u \Ge \theta v.$
\end{enumerate}

Let $u$ be a derivative,
that is, $u = \theta y_j$ for a derivative operator 
$$\theta = \delta_1^{k_1}\delta_2^{k_2}\cdots\delta_m^{k_m} \in \Theta$$
 and $1\Le j \Le n$. The {\it order} of $u$ is defined as 
$$\ord u=\ord\theta=k_1+\ldots+k_m.$$ If $f$ is a differential 
polynomial, $f\not\in\K$, then $\ord f$ denotes the maximal order of 
derivatives appearing effectively in $f$. 

A ranking $>$ is called
{\it orderly} iff $\ord u > \ord v$ implies $u > v$ for all 
derivatives $u$ and $v$. A ranking $>_{el}$ is called 
an {\it elimination} ranking iff $y_i >_{el} y_j$ implies 
$\theta_1y_i >_{el} \theta_2y_j$ for all 
$\theta_1, \theta_2 \in \Theta$.

Let a ranking $<$ be fixed.
The derivative $\theta y_j$ of the highest rank appearing 
in a differential polynomial $f \in \K\{y_1,\dots,y_n\} \setminus \K$ 
is called the {\it leader} of $f$. We denote the leader by $\ld f$ or $\u_f$.
The indeterminate $y_j$ is called the {\it leading variable} of $f$ and
denoted by $\lv f.$ 
Represent $f$ as a univariate polynomial in $\u_f$:
$$
f = \i_f \u_f^d + a_1 \u_f^{d-1} + \ldots + a_d.
$$
The monomial $\u_f^d$ is called the {\it rank} of $f$ and
is denoted by $\rk f.$
Extend the ranking relation on derivatives variables to ranks: 
$u_1^{d_1}>u_2^{d_2}$ iff either $u_1>u_2$ or
$u_1=u_2$ and $d_1>d_2$.

The polynomial $\i_f$ is called the {\it initial} of $f$.
Apply any $\delta \in \Delta$ to $f$:
$$ \delta f = \frac{\partial f}{\partial \u_f}\delta \u_f + \delta \i_f \u_f^d +
\delta a_1 \u_f^{d-1}+\ldots + \delta a_d.
$$
The leader of $\delta f$ is $\delta \u_f$ and the initial 
of $\delta f$ is called the {\it separant} of $f$, denoted $\s_f$. 
If $\theta\in\Theta\setminus\{1\}$, then $\theta f$ is called a
{\it proper derivative} of $f$. Note that the initial of any 
proper derivative of $f$ is equal to $\s_f$. 

We say that a differential polynomial $f$ is {\it partially reduced} 
w.r.t. $g$ iff no proper derivative of $\u_g$ appears in $f$. 
A differential polynomial $f$ is {\it algebraically reduced} w.r.t.
$g$ iff $\deg_{\u_g} f<\deg_{\u_g}g$. 
A differential polynomial $f$ is {\it reduced} w.r.t. a differential
 polynomial $g$ iff $f$ is partially and algebraically reduced
 w.r.t. $g$. Consider any subset 
$\mathbb{A} \subset \K\{y_1,\ldots,y_n\}\setminus \K$. We say that 
$\mathbb{A}$ is autoreduced (respectively, algebraically autoreduced) 
iff each element of $\mathbb{A}$ is reduced (respectively, algebraically reduced) 
w.r.t. all the others. 

Every autoreduced set is finite 
\cite[Chapter I, Section 9]{Kol} (but an algebraically autoreduced set
in a ring of differential polynomials
may be infinite). For autoreduced sets we use capital letters 
$\mathbb{A, B, C,}$ \ldots and 
notation $\mathbb{A} = A_1,\ldots,A_p$ to specify the list of the elements 
of $\mathbb{A}$ arranged in order of increasing rank. 

We denote the sets of initials and separants of 
elements of $\mathbb{A}$ by $\i_\mathbb{A}$ and $\s_\mathbb{A}$,
 respectively. Let $H_A=\i_\mathbb{A}\cup \s_\mathbb{A}$.
Let $S$ be a finite set of differential polynomials. 
Denote by $S^\infty$ the multiplicative set containing $1$ and
 generated by $S$. Let $I$ be an ideal in a commutative ring $R$. 
The {\it saturated ideal} $I:S^\infty$ is defined as 
$\{a \in R\:|\:\exists s \in S^\infty: sa \in I\}$. If $I$ is a 
differential ideal then $I:S^\infty$ is also a differential ideal 
(see \cite{Kol}).

Consider two polynomials $f$ and $g$ in
 $\K\{y_1,\ldots, y_n\}$. Let $I$ be the differential ideal 
generated by $g$. Applying a finite number of pseudo-divisions,
one can compute a {\it differential partial remainder} $f_1$ and a 
{\it differential remainder} $f_2$
of $f$ w.r.t. $g$ such that there exist $s \in S_g^\infty$ and
 $h \in H_g^\infty$ satisfying $sf \equiv f_1$ and 
$hf \equiv f_2 \mod I$ with $f_1$ and $f_2$ partially reduced 
and reduced w.r.t. $g$, respectively (see \cite{Fac} for 
definitions and the algorithm for computing remainders). We denote by
$\drem(f,\A)$ the differential remainder of a polynomial $f$
w.r.t. an autoreduced set $\A$.

Let $\mathbb{A} = A_1,\ldots,A_r$ and $\mathbb{B} = B_1,\ldots,B_s$
 be (algebraically) autoreduced sets. We say that $\mathbb{A}$ has lower rank than 
$\mathbb{B}$ if
\begin{itemize}
\item
there exists $k \Le r, s$ such that $\rank A_i$ = 
$\rank B_i$ for $1 \Le i < k,$ and $\rank A_k < \rank B_k$,
\item or if $r > s$ and $\rank A_i = \rank B_i$ for $1 \Le i \Le s$. 
\end{itemize}
We say that
 $\rank\mathbb{A} = \rank\mathbb{B}$ iff $r=s$ and 
$\rank A_i = \rank B_i$ for $1 \Le i \Le r$.

The following notion of a characteristic set in {\it Kolchin's sense} 
in characteristic zero is crucial in our further discussions.
 It was first introduced by Ritt for prime differential ideals, and 
then extended by Kolchin to arbitrary differential ideals.

\begin{defn}\cite[page 82]{Kol} An autoreduced subset of the
  lowest rank in a set $X\subset \K\{Y\}$ is
called a {\it characteristic set} of $X$.
\end{defn}

We call these sets {\it Kolchin characteristic sets} to avoid 
confusion with other notions, e.g., in \cite{Fac,Dif} characteristic 
sets are used in Kolchin's sense and in some other senses. 
A characteristic set in Kolchin's sense exists for any 
set $X\subset\K\{Y\}$ due to the fact that every family
of autoreduced sets contains one of the least rank (see \cite{Kol}).

As it is mentioned in \cite[Lemma 8, page 82]{Kol}, in the case 
of $\Char k = 0$, a set $\mathbb{A}$ is a characteristic set of a 
proper differential ideal $I$ iff each element of $I$ reduces to 
zero w.r.t. $\mathbb{A}$. Moreover, the leaders and the correspondent
degrees of these leaders of any two characteristic sets of $I$
coincide.

\begin{defn}\cite[Definition 2.6]{Fac} A differential ideal $I$ in $\K\{y_1,\ldots,y_n\}$ is
said to be {\it characterizable} if there exists a characteristic set $\mathbb{A}$ of $I$ in Kolchin's sense
such that $I = [\mathbb{A}]:H_\mathbb{A}^\infty.$ We call any such characteristic set $\mathbb{A}$ a {\it characterizing} set
of $I$.
\end{defn}

Characterizable ideals are radical \cite[Theorem 4.4]{Fac}.

\section{Rosenfeld-Gr\"obner algorithm for the ordinary case} \label{sec:rg}

A system of ordinary differential equations and inequalities $\A=0,H\neq 0$,
where $\A,H\subset\K\{Y\}$,
is called regular (see \cite{Bou1}),
if $A$ is autoreduced, $H$ is partially reduced w.r.t. $\A$, and 
$H\supseteq H_\A$, where $H_\A$ is the set of initials and separants
of elements of $\A$ (in the partial differential case it is also
required that the set $\A$ is coherent, but in the ordinary case
this condition holds for any autoreduced set $\A$).
For a regular system $\A,H$, the differential ideal
$[\A]:H^\infty$ is also called regular. Every regular ideal is radical (see \cite{Bou1}), 
and, according to the Rosenfeld Lemma, $f\in [\A]:H^\infty$ if and only
if the partial remainder of $f$ w.r.t. $\A$ belongs to the algebraic
ideal $(\A):H^\infty$.

The Rosenfeld-\G\ algorithm proposed in \cite{Bou1,Bou2} computes a
regular decomposition of a given radical differential ideal $\{F\}$, i.e.,
a representation $$\{F\}=\bigcap_{i=1}^k[\A_i]:H_i^\infty,$$ where
$[\A_i]:H_i^\infty$ are regular differential ideals.

 We begin with the following version of the Rosenfeld-Gr\"obner 
 algorithm. It is very similar to the original algorithm
 presented in \cite{Bou1}, except for 
 the fact that we are in the ordinary case and need not
 deal with coherence. We also note that some of the regular systems 
 computed by the version of the algorithm presented here may
 correspond to unit ideals; this can be checked later on by 
 means of Gr\"obner basis computations as in \cite{Bou1}
 or via polynomial GCD computations modulo regular chains
 as in \cite{Bou3}.
 
Finally, we follow the suggestion 
 given in 
 \cite[Improvements, page 73]{Dif}: it is recommended to reduce 
 the multiplicative set $H$ of initials and separants. If it turns 
 out that one of them reduces to zero, then the corresponding saturated component
 contains $1$ and therefore need not be considered. We implement these ideas
in Algorithm~\ref{RGconv}.

\begin{figure}[htb]
\begin{alg}\label{RGconv}{\sf Rosenfeld-Gr\"obner}$(F_0,H_0)$\\
\begin{tabular}{ll}
{\sc Input:}& finite sets of differential polynomials $F_0,H_0$ \\
            & and a differential ranking\\
{\sc Output:} & a finite set $T$ of regular systems such that \\
& $\qquad\qquad\{F_0\}:H_0^\infty=\bigcap\limits_{(\A,H)\in
  T}[\A]:H^\infty\;$ 
\end{tabular}\\
\begin{tabular}{l}
\b $T:=\varnothing$, $\;\;\;U:=\{(F_0,H_0)\}$\\
\b {\bf while} $U\neq\varnothing$ {\bf do}\\
\b \b Take and remove any $(F,H)\in U$\\
\b \b $\C:=$ characteristic set of $F$\\
\b \b $\bar F:=\Rem(F\setminus\C,\C)\setminus\{0\}$\\
\b \b $\bar H:=\Rem(H,\C)\cup H_\C$\\
\b \b {\bf if} $\bar F\cap\K=\varnothing$ {\bf and} $0\not\in\bar H$ {\bf then}\\
\b \b\b {\bf if} $\bar F=\varnothing$ {\bf then} $T:=T\,\cup\,\{(\C,\bar H)\}$\\
\b \b\b\b {\bf else} $U:=U\cup\{(\bar F\cup\C,\bar H)\}$\\
\b \b\b {\bf end if}\\
\b \b {\bf end if}\\
\b \b $U:=U \cup\{(F\cup\{h\},H)\;|\;h\in H_\C,\;h\not\in\K\cup H\}$\\
\b {\bf end while}\\
\b {\bf return} $T$
\end{tabular}
\end{alg}
\end{figure}

 Given a set $F$ of differential polynomials,
 the Rosenfeld-Gr\"obner algorithm at first computes a characteristic set $\C$ of
 $F$, i.e., an  autoreduced subset of $F$ of the least rank. 
 It may happen that $\lv\C\subsetneq\lv F$ (for example, take
 $F=\{x+y,y\}$ w.r.t. a ranking such that $x>y$). In other words,
 inclusion $F_1\subset F_2$ does not imply that for the corresponding
 characteristic sets $\C_1$ and $\C_2$, we have $\C_1\subseteq\C_2$.
 We need the latter property, in order to obtain the bound, 
 so we are going to relax the requirement that $\C$ is
 autoreduced.

 A subset $\C$ of $\K\{Y\}\setminus\K$ is called a weak d-triangular 
 set \cite[Definition 3.7]{Dif}, if the set of its leaders $\ld \C$ is 
 autoreduced. In the ordinary case, $\C$ is a weak d-triangular set if 
 and only if the leading differential indeterminates $\lv f$, $f\in\C$, 
 are all distinct. A partially autoreduced weak d-triangular set is called
 d-triangular \cite[Definition 3.7]{Dif}.
 For a polynomial $f$ and a weak d-triangular
 set $\C$, the pseudo-remainder $\Rem(f,\C)$ is defined
 via \cite[Algorithm 3.13]{Dif}.

We will replace the reduction of $F$ w.r.t. an autoreduced set
in the Rosenfeld-Gr\"obner algorithm by that w.r.t. a weak
d-triangular set. We note that the version of the algorithm
presented in \cite[Section 6]{Dif} (Algorithms 6.8, 6.10, and 6.11)
also computes differential pseudo-remainders w.r.t.  weak
d-triangular sets. Since the output regular systems must be 
partially autoreduced, at the very end, partial autoreduction of
the weak d-triangular set $\C$ via
\cite[Algorithm 6.8]{Dif} is carried out. 

Alternatively,
one could perform partial autoreduction every time 
a weak d-triangular set is updated. In the following
section, we show how to perform this autoreduction,
as well as computation of differential pseudo-remainders, 
so that the inequality $$M(F\cup H)\Le (n-1)!M(F_0\cup H_0)$$ is preserved
(see formula \eqref{BoundDefinition} below).

\section{Modified Rosenfeld-\G\ algorithm}\label{sec:modifiedRG}

For a set of differential polynomials $F$, let 
$m_i(F)$ be the maximal order of the differential
indeterminate $y_i\in Y$ occurring in $F$. If $y_i$ does not
occur in $F$, we set $m_i(F)=0$.
Let 
\begin{align}\label{BoundDefinition}
M(F)=\sum_{i=1}^n m_i(F).
\end{align}
We propose a modification of the Rosenfeld-\G\ algorithm 
(see Algorithm~\ref{RGBound} below), in
which for every intermediate system $(F,\C,H)\in U$, the bound
\begin{align}\label{formula:bound}
M(F\cup\C\cup H)\Le (n-1)!M(F_0\cup H_0)
\end{align} 
holds, where $F_0=0,H_0\neq 0$ is
the input system of equations and inequalities
corresponding to the radical differential ideal $\{F_0\}:H_0^\infty$.

In the formula~\eqref{formula:bound} we have a multiple $(n-1)!.$
If the number of variables is equal to $1$ or $2$ it disappears.
In the case of $n=2$ Ritt proved the Jacobi bound for 
$|F_0| = 2$ and empty $H_0$ by the direct computation and his
result does not have any multiple either. Consider the intuition 
behind the case of $n = 3$ by looking at a particular example.

\begin{exmp} Let $F_0 = x+y+z,$ $x'$ with
the elimination ranking $x>y>z.$ Then $m_x = 1,$ $m_y = m_z = 0$
and $$M(F_0) = 1+0+0 = 1.$$ 
In order to find a characteristic set of the prime differential
ideal $[F_0]$ we reduce $x'$ w.r.t. $x+y+z$ and get $y'+z'.$
The output consists of two polynomials: $$\C = y'+z', x+y+z.$$
We have: $m_x(\C) = 0,$ $m_y(\C) = 1,$ and $m_z(\C) = 1.$
Hence, $$M(\C) = 0+1+1 =2 > 1.$$ But $(n-1)! = (3-1)! = 2! = 2$ and
$2 \Le 2\cdot 1.$
\end{exmp}

\subsection{Algebraic computation of differential remainders}\label{sec:da}

The Rosenfeld-Gr\"obner algorithm requires to compute differential
pseudo-remainders $R=\Rem(F\setminus\C,\C)$.
If the ranking on derivatives is not orderly, the orders
of some (non-leading) derivatives
may grow as a result of the differential pseudo-reduction, so that we may
have $m_i(R)>m_i(F)$ for some $i\in\{1,\ldots,n\}$. To ensure a 
bound on $m_i(R)$, we construct a triangular 
set\footnote{A set is called triangular if the leaders of its
elements are distinct.} $\B$,
such that the computation of the differential pseudo-remainders $\Rem(F,\C)$ 
can be replaced by the computation of algebraic pseudo-remainders $\algrem(F,\B)$, 
and, at the same time, $\B$ satisfies a bound on the orders of 
derivatives occurring in it. 

For a set $\B$ of differential polynomials and a differential indeterminate
$v\in\lv\B$, let 
$$
\B_v=\{f\in\B\;|\;\lv f=v\}.
$$
Assume that $\B$ is algebraically triangular, which implies that 
for any non-empty subset $\A\subset\B$, elements
of $\A$ of the minimal and maximal ranks are uniquely defined
and denoted, respectively, $\min\A$ and $\max\A$.
Define the following two subsets of $\B$:
\begin{align*}
\B^0 &= \{\min\B_v\;|\;v\in\lv\B\}\\
\B^* &= \{\max\B_v\;|\;v\in\lv\B\}.
\end{align*}

Also, for a set $\{m_i\}_{i=1}^k$ of non-negative integer numbers
and an arbitrary set $F$ of differential polynomials, let
$$
F_{\{m_i\}}=\{f\in F\;|\;\ord_{y_i}f\Le m_i, \;i=1,\ldots,k\}.
$$

\begin{figure}
\begin{alg}\label{DifferentiateAutoreduce}{\sf Differentiate\&Autoreduce}$(\C,\{m_i\})$\\
\begin{tabular}{l}
{\sc Input:} a weak d-triangular set $\C=C_1,\ldots,C_k$ with 
              $\ld \C=y_1^{(d_1)},\ldots,y_k^{(d_k)}$,\\
\hphantom{\sc Input:} and a set of non-negative integers $\{m_i\}_{i=1}^k$, $m_i\Ge m_i(\C)$ \\
{\sc Output:}  set $\B=\left\{B_i^j\;\big|\;1\Le i\Le k,\;0\Le j\Le m_i-d_i\right\}$ satisfying\\
\hphantom{\sc Output:} \b\b $\rk B_i^j=\rk C_i^{(j)}$ \\
\hphantom{\sc Output:} \b\b$B_i^j$ are 
reduced w.r.t. $\C\setminus\{C_i\}$\\
\hphantom{\sc Output:} \b\b$m_i(\B)\Le m_i,$ $i = 1,\ldots,k$\\
\hphantom{\sc Output:} \b\b$m_i(\B)\Le m_i(\C)+\sum_{j=1}^k (m_j-d_j)$, $i=k+1,\ldots,n$ \\
\hphantom{\sc Output:} \b\b $\B\subset\left[\B^0\right]\subset[\C]\subset[\B]:H_\B^\infty$\\
\hphantom{\sc Output:} \b\b $H_\B\subset H_\C^\infty+[\C],$
                  $\;\;H_\C\subset (H_\B^\infty+[\B]):H_\B^\infty$\\
\hphantom{\sc Output:} {\bf or} $\{1\}$, if it is detected that $[\C]:H_\C^\infty=(1)$
\end{tabular}\\
\begin{tabular}{rl}
1&  $\D:=\C$, $\B:=\varnothing$\\
2&  {\bf while} $\D\cup(\delta\B^*)_{\{m_i\}}\neq\varnothing$ {\bf do}\\
3&  \b $f:=\min \left(\D\cup(\delta\B^*)_{\{m_i\}}\right)$\\
4&  \b {\bf if} $f\in\D$ {\bf then} \\
5&  \b\b $\bar f:=\algrem(f,\B)$\\
6&  \b\b $\D:=\D\setminus\{f\}$\\
7&  \b {\bf else} \\
8&  \b\b$\bar f:=\algrem(f,\B\cup(\delta\B^0\setminus\{f\}))$\\
9&  \b {\bf end if}\\
10&  \b {\bf if} $\rk\bar f\neq\rk f$ {\bf then return} $\{1\}$
       {\bf end if}\\
11&  \b $\B:=\B\cup\left\{\bar f\right\}$\\
12& {\bf end while} \\
13& {\bf return} $\B$
\end{tabular}
\end{alg}
\end{figure}

Before we prove correctness and termination of 
Algorithm {\sf Dif\-fe\-ren\-ti\-ate\&Au\-to\-re\-du\-ce}, let us
discuss it informally. 
The triangular set $\B$
computed by the algorithm can be thought of as a result of an autoreduction 
of a {\it differential prolongation} of the input set $\C=\{C_1,\ldots,C_k\}$, i.e., 
of the set $$\tilde\C=\{\delta^jC_i\;|\;1\Le i\Le k,\;0\Le j\Le
m_i-d_i\}.$$
In particular, we have $\rk\B=\rk\tilde\C$, unless the autoreduction
process cancels one of the initials, in which case we can show
that $[\C]:H_\C^\infty=(1)$. 

However, if one wants to make 
this autoreduction completely algebraic
(in order to control the growth of orders), one has to be careful, 
because in the above set $\tilde\C$ there may appear derivatives
of some $\ld C_i$ of order higher than those that appear in $\ld\tilde\C$,
which cannot be canceled 
by an algebraic reduction. For example,
if $\C=\{y_1,y_2+y_1'\}$, $m_1=1$, $m_2=2$, and the ranking is
elimination with $y_1<y_2$, then 
$$\tilde\C=\{y_1,y_1',y_2+y_1',y_2'+y_1'',y_2''+y_1'''\},$$
and in the last two polynomials derivatives $y_1'',y_1'''$ cannot
be canceled by algebraic reduction w.r.t. $y_1$ and $y_1'$.

This problem is avoided by computing the elements of $\B$ in the 
order of increasing rank. If the polynomials are added to $\B$ in this
order, one only needs to reduce each new polynomial $f$, which 
we are going to add to $\B$, w.r.t. the set $\delta\B^0$ of first-order 
derivatives of the elements of $\B^0$ (this set $\B^0$ has the same leaders as $\C$) and the set $\B$; this will guarantee that $f$ 
is differentially  
reduced w.r.t. $\B^0\setminus \B_{\lv f}.$ 

The inclusions
$\B\subset[\B^0]\subset[\C]\subset[\B]:H_\B^\infty$,
$H_\B\subset H_\C^\infty+[\C]$, and $H_\C\subset (H_\B^\infty+[\B]):H_\B^\infty$
will allow us to replace reduction w.r.t. $\C$ (in the 
Rosenfeld-Gr\"obner algorithm) by that 
w.r.t. $\B$ without disturbing the saturated ideal. 
If $m_i$'s are chosen as the maximal orders of derivatives of $y_i$'s,
$1\Le i\Le k$,  appearing in the set that is being reduced w.r.t. $\C$, 
then we can replace the differential reduction w.r.t. $\C$ by
the algebraic reduction w.r.t. $\B$. The orders of derivatives 
of $y_i$'s appearing in the remainder then will not exceed
$d_i$ for the leading $y_i$'s (i.e., for $1\Le i\Le k$). For the 
non-leading $y_i$'s, the orders are bounded by the inequality
\begin{align}\label{BInequality} 
m_i(\B)\Le m_i(\C)+\sum_{j=1}^k (m_j-d_j),\;\;\;i=k+1,\ldots,n.
\end{align}

We will use the following two lemmas in the proof of correctness
of Algorithm {\sf Differentiate\&Autoreduce}. 

\begin{lem} \label{InitialBelongs} 
  Let $\C$ be a weak d-triangular set in the ring of differential
  polynomials $\K\{Y\}$ with derivations
  $\Delta=\{\delta_1,\ldots,\delta_m\}$. Assume that a ranking
  on the set of derivatives $\Theta Y$ is fixed. 
  Let $f\in\K\{Y\}$ be a differential
  polynomial with $\ld f\not\in\Theta\ld \C$, and let $f\to_\C g$. Then
  \begin{itemize}
     \item $\rk g<\rk f\;\Rightarrow\;\i_f\in [\C]:H_\C^\infty$
     \item $\rk g=\rk f\;\Rightarrow\;\exists\;h\in H_\C^\infty\;$
     such that $h\cdot\i_f-\i_g\in [\C]$, $h\cdot\s_f-\s_g\in [\C]$.
  \end{itemize}
\end{lem}
\begin{pf}
  Let $\rk f=u^d$, and let $\A=\{p\in\Theta\C\;|\;\ld p<u\}$. Then for
  every $p\in\A$, $p$ and $\i_p$ are free of $u$.

  Since $f\to_\C g$ and $u\not\in\Theta\ld\C$, 
there exist polynomials
  $h\in\i_\A^\infty$, $A_1,\ldots,A_k\in\A$ and
  $\alpha_a,\ldots,\alpha_k\in\K\{Y\}$
  such that
  \begin{equation} \label{eqred}
    h\cdot f=g+\sum_{i=1}^k\alpha_iA_i.
  \end{equation}
  The maximal degree of $u$ present in \eqref{eqred} is equal to
  $d$. {\it Replace} every occurrence of $u^d$ by a new variable $v$, 
  and consider \eqref{eqred} as an equality between two polynomials
  in $v$, in which polynomials $h,A_1,\ldots,A_k$ are free of $v$.
  We have therefore:
  $$h\cdot\frac{df}{dv}=\frac{dg}{dv}+\sum_{i=1}^k\frac{d\alpha_i}{dv}A_i.$$
  It remains to notice that $\frac{df}{dv}=\i_f$ and
  $$\frac{dg}{dv}=\left\{\begin{array}{ll}
       0,&\rk g<\rk f\\
       \i_g,& \rk g=\rk f\end{array}\right.,
  $$
  hence we obtain
  \begin{itemize}
    \item $\rk g<\rk
    f\;\Rightarrow\;\i_f\in(\A):\i_\A^\infty\subset[\C]:H_\C^\infty$.
    \item $\rk g=\rk
    f\;\Rightarrow\;h\cdot\i_f-\i_g\in(\A)\subset[\C]$, where 
    $h\in\i_\A^\infty\subset H_\C^\infty$.
  \end{itemize}
 
  Consider now \eqref{eqred} as an equality between two polynomials in
  $u$, in which $h,A_1,\ldots,A_k$ are free of $u$. We have therefore:
  $$h\cdot\frac{df}{du}=\frac{dg}{du}+\sum_{i=1}^k\frac{d\alpha_i}{du}A_i.$$
  It remains to notice that $\frac{df}{du}=\s_f$ and, if $\rk g=\rk
  f$, $\frac{dg}{du}=\s_g$, hence $h\cdot\s_f-\s_g\in(\A)\subset[\C]$,
  where $h\in\i_\A^\infty\subset H_\C^\infty$.
\end{pf}

\begin{rem}\label{InitialRemark} The above lemma also holds when the set of
  derivations $\Delta$ is empty, in which case $\K\{Y\}=\K[Y]$ is a
  ring of algebraic polynomials, $\C\subset\K[Y]$ is a triangular set,
  $\to_\C$ is the algebraic pseudo-reduction relation w.r.t. $\C$, 
  and $[\C]=(\C)$ is an ideal in $\K[Y]$.
\end{rem}

\begin{lem} \label{l:Hubert} \cite[Lemma 6.9]{Dif}
  Let $H$ and $K$ be two sets of differential polynomials, and 
  let $I$ be a differential ideal. If $K\subset
  (H^\infty+I):H^\infty$, then $I:H^\infty=I:(H\cup K)^\infty$.
\end{lem}
\begin{pf}
  The proof of this statement is omitted in \cite{Dif}, so,
  for the sake of completeness, we provide it here.

  Clearly, $I:H^\infty\subseteq I:(H\cup K)^\infty$. To prove the
  inverse inclusion, take any $f\in I:(H\cup K)^\infty$. 
  Then, by definition, there exist $h\in H^\infty$ and $k\in K^\infty$
  such that $fhk\in I$. 

  Since $K\subset  (H^\infty+I):H^\infty$, there exist $h_1,h_2\in
  H^\infty$ such that $kh_1-h_2\in I$. The fact that $fhk\in I$ 
  implies that $fhkh_1\in I$, whence $$fhkh_1-fh(kh_1-h_2)=fhh_2\in I,$$ i.e.,
  $f\in I:H^\infty$. 
\end{pf}

\begin{prop} Algorithm {\sf Differentiate\&Autoreduce} is correct and terminates.
\end{prop}
\begin{pf}
The proof of specifications of the algorithm is divided into 
{\bf three} parts. We will:
\begin{itemize}
\item first prove the statements about the ranks
of the elements of $\B$, 
\item then about their orders, and 
\item finally, the inclusions. 
\end{itemize}
All these statements hold {\it only} if the condition in line 10 is never satisfied, that is, throughout 
the algorithm $\rk\bar f=\rk f$; for now we assume this. At 
the end of this proof we will show that, if $\rk\bar f\neq\rk f$
for some $f$, then the ideal $[\C]:H_\C^\infty$ must be trivial.
Let us proceed to the three parts of the proof:
\begin{enumerate}
\item[\bf 1.] We prove that the output $\B$ has the 
form 
$$\B=\left\{B_i^j\;\big|\;1\Le i\Le k,\;0\Le j\Le m_i-d_i\right\},$$
where $\rk B_i^j=\rk C_i^{(j)}$. 
For $i=1,\ldots,k$, let
$$j_i=\left\{\begin{array}{ll}
     m_i(\ld\B)-d_i,& y_i\in\lv\B\\
     -1,&{\rm otherwise}.
    \end{array}\right.$$

The statement of Part 1 follows from the following invariants of 
the {\bf while}-loop:
\begin{align}
&-1\Le j_i\Le m_i-d_i,\ i=1,\ldots,k\tag{I1}\label{Iji}\\
&\B=\left\{B_i^j\;\big|\;1\Le i\Le k,\;0\Le j\Le j_i\right\}\tag{I2}\label{IB}\\
&\rk B_i^j=\rk C_i^{(j)},\  \left(1\Le i\Le k,\;0\Le j\Le j_i\right)\tag{I3}\label{IrkB} \\
&\lv\B\cap\lv\D=\varnothing\tag{I4}\label{Ilv}\\
&\D\subset\C\tag{I5}\label{IDC}\\
&\text{For all}\ f\in\B, g\in\D, \rk f<\rk g.\tag{I6}\label{BlessD}
\end{align}

One can check immediately that the above invariants hold at the beginning 
of the first iteration 
of the {\bf while}-loop. Assume that we are at the beginning of some 
iteration and the invariants hold; show that they will also hold at the end of 
this iteration.

Let $f$ be the polynomial computed in line 3, and let $y_i^{(d)}=\ld f$. 
We have {\bf two cases}:
\begin{itemize}
\item $f\in\D$. By \ref{IDC} and the fact that the 
leading variables of the elements of $\C$ are distinct, we have 
$\rk f=\rk C_i$. By \ref{Ilv} $\lv f\not\in\lv\B$, whence by definition 
of $j_i$ we have $j_i=-1$. Since $\rk\bar f=\rk f$, at the end of the
iteration we will have $B_i^0=\bar f$ with $\rk B_i^0=\rk C_i$, and
$j_i=0$. Thus, invariants \ref{Iji}--\ref{IrkB} will hold. 

Invariants \ref{Ilv} and \ref{IDC} also continue to hold due to 
the assignments in lines 6 and 11 and the fact that sets $\D$ and
$\B$ do not change elsewhere throughout the iteration of the {\bf while}-loop.
The choice of $f$ in line 3 and 
the assignment in line 6 also imply that at the end of the iteration
we have $\rk\bar f=\rk f<\rk g$ for all $g\in\D$, whence invariant
\ref{BlessD} is preserved. 
\item  $f\in(\delta\B^*)_{\{m_i\}}$. By
\ref{IB} and \ref{IrkB}, $j_i\Ge 0$ and $\rk f=\rk C_i^{(j_i+1)}$. 
Hence, at the end of the iteration $j_i$ increases by one, 
while $\bar f$ with $\rk\bar f=\rk f$ is added to $\B$,
thus preserving invariants \ref{Iji}--\ref{IrkB}. Note
also that $\lv f\in\lv\B$, whence $\lv\B$ is preserved as well.
Hence, due to the fact that $\D$ remains unchanged throughout
the iteration, invariants \ref{Ilv} and \ref{IDC} are preserved.
The facts that $\rk\bar f=\rk f\le \rk g$ for all $g\in\D$
(due to the choice of $f$ in line 3)
and that $\lv\bar f\not\in\lv\D$ (due to \ref{Ilv}) implies
preservation of \ref{BlessD} at the end of the iteration.
\end{itemize}

The above also proves the {\bf termination} of the algorithm: at each 
iteration exactly one of the $j_i$ is incremented, whence 
the number of iterations does not exceed 
$$\sum_{i=1}^k (m_i-d_i+1).$$

\item[\bf 2.] We first show that:
\begin{itemize}
\item the elements $B_i^j$ of the 
output $\B$ are reduced w.r.t. $\C\setminus\{C_i\}$ and that
\item 
$m_i(\B)\Le m_i$, $i=1,\ldots,k$. 
\end{itemize}
This is implied by the fact that $m_i(\rk\B)\Le m_i$, 
$i=1,\ldots,k$, which is a consequence of Part 1, and the 
following invariants:

\begin{align} 
&m_t\left(B_i^j\right)\Le d_t,\ (1\Le i\neq t\Le k,\; 0\Le j\Le j_i)\tag{I7}\label{IBmt}\\
&B_i^j\ \text{are differentially reduced w.r.t.}\ \B^0\setminus\{B_i^0\} \tag{I8}\label{IBred}
\end{align}

We have {\bf two cases}:
\begin{itemize}
\item $f\in\D$. Let us show that
\begin{equation} \label{eq:mtf}
m_t(f)\Le d_t+j_t,\quad y_t\in\lv\B.
\end{equation}
The fact that $y_t\in\lv\B$, according to \ref{IrkB}, 
implies $j_t\Ge 0$. 
We may assume that $y_t$ is present in 
$f$: otherwise $m_t(f)=0$ and (\ref{eq:mtf}) will
trivially hold due to the fact that $d_t\Ge 0$ and $j_t\Ge 0$.
According to \ref{Iji}, two cases are 
possible:

\begin{enumerate}
\item $0\Le j_t<m_t-d_t$. Then there exists a polynomial
$g\in(\delta\B)^*_{\{m_i\}}$ with $\lv g=y_t$. By \ref{IrkB}, 
$\ld g=y_t^{(d_t+j_t+1)}$. Since, due to \ref{Ilv}, $y_t$ cannot be
the leading variable of $f$, yet $y_t$ is present in $f$,
$y_t^{(m_t(f))}<\ld f$. Since $f$ is an
element of $\D\cup(\delta\B)^*_{\{m_i\}}$ of the least rank, $\ld f\le \ld g$. 
Combining these statements, we obtain
$$
y_t^{(m_t(f))}<\ld f\le \ld g=y_t^{(d_t+j_t+1)},
$$
which implies (\ref{eq:mtf}).
\item $j_t=m_t-d_t$. Then, due to \ref{IDC}
and the condition on the input $\C$ we have 
$m_t(f)\Le m_t$, which yields (\ref{eq:mtf}).
\end{enumerate}

Inequality (\ref{eq:mtf}) and invariant \ref{IBmt} imply that 
the algebraic remainder $\bar f$ computed in line 5 is
differentially reduced w.r.t. $\B^0$ and satisfies
\begin{equation} \label{eq:mtbarf}
m_t\left(\bar f\right)\Le d_t,\quad y_t\in\lv\B.
\end{equation}
Note also that due to \ref{BlessD}, $\B$ is differentially 
reduced w.r.t. $\bar f$. Thus, 
invariant \ref{IBred} also holds at the end of the 
iteration.

Taking into account that for all $g\in\D$ we 
have $\rk\bar f=\rk f\Le\rk g$, we obtain
\begin{equation} \label{eq:mtbarf1}
  m_t\left(\bar f\right)\Le d_t,\quad y_t\in\lv\D.
\end{equation}
Together inequalities (\ref{eq:mtbarf},\ref{eq:mtbarf1}) yield
invariant \ref{IBmt} at the end of the iteration.

\item $f\in(\delta\B^*)_{\{m_i\}}$. 
By \ref{IBmt}, $m_t(f)\Le d_t+1$, for $t$ such that 
$y_t\in\lv\B\setminus\{\lv f\}$.
This inequality and invariant \ref{IBmt} imply that the algebraic remainder 
$\bar f$ computed in line 8 is differentially reduced w.r.t. 
$\B^0\setminus\{B_i^0\}$ and satisfies $m_t(\bar f)\Le d_t$.
Thus, we obtain that invariants \ref{IBmt} and 
\ref{IBred} also holds at the end of the iteration. 
\end{itemize}

Finally, we prove the {\bf bound for the orders} of non-leading
derivatives in the output:
$$m_s(\B)\Le m_s(\C)+\sum_{t=1}^k (m_t-d_t),\;\;s=k+1,\ldots,n.$$
This bound holds due to the following invariant:
\begin{equation} \label{invbound}
m_s(\B)\Le m_s(\C)+\sum_{\{t: j_t\Ge 0\}} j_t,\;\;s=k+1,\ldots,n.
\end{equation}

Assume that (\ref{invbound}) holds at the beginning of an iteration, 
and let $s\in\{k+1,\ldots,n\}$. We then have:
\begin{enumerate}
\item
If $f\in\D$, no differentiations occur during the iteration
and the sum remains unchanged, whence (\ref{invbound}) is preserved.
\item
If $f\in(\delta\B^*)_{\{m_i\}}$, then $f=\delta g$ for some $g\in\B$,
whence $m_s(f)\Le m_s(\B)+1$. Similarly, 
$$
m_s\left(\B\cup\delta\B^0\right)\Le m_s(\B)+1.
$$ 
Thus, according to line 8, $m_s\left(\bar f\right)\Le m_s(\B)+1$, 
and so at the end of the iteration $m_s(\B)$ is increased
at most by one. At the same time, as was shown in Part 1, 
Case 2, exactly one of the $j_i$ is incremented, whereby
the sum in (\ref{invbound}) increases by 1. Thus, (\ref{invbound})
is preserved. 
\end{enumerate}

\item[\bf 3.] It remains to prove:
\begin{enumerate}
\item the inclusions in
the specification of the algorithm and that
\item whenever the algorithm
outputs $\{1\}$, the ideal $[\C]:H_\C^\infty$ is trivial.
\end{enumerate}

The inclusions are implied by the following invariants:
\begin{align}
&\B\subset\left[\B^0\right]\subset [\C]\tag{I9} \label{IB0}\\
&H_\B\subset H_\C^\infty+[\C]\tag{I10} \label{IHB}\\
&\C\setminus\D\subset [\B]:H_\B^\infty\tag{I11} \label{IC}\\
&H_{\C\setminus\D}\subset (H_\B^\infty+[\B]):H_\B^\infty\tag{I12} \label{IHC}
\end{align}

Assume that the invariants hold at the beginning of some iteration; show
that either the algorithm terminates at this iteration with the output
$\{1\}$ or the invariants will hold at the end of this iteration.
We have {\bf two cases}:

\begin{itemize}
\item $f\in\D$. Then by \ref{IDC} $f\in\C$. 
Since $\bar f=\algrem(f,\B)$, we have 
$\bar f\in (\B\cup\{f\})$. Then, according to \ref{IB0}, $\bar f\in[\C]$.
As was shown in Part 1, Case 1, $\bar f$ is added to $\B^0$ in line 11,
thus preserving \ref{IB0}.

Next, due to \ref{Ilv}, $\ld f\not\in\ld\B$. Thus, 
Lemma~\ref{InitialBelongs} (see also Remark~\ref{InitialRemark}) 
applies to the algebraic remainder 
computed in line 5. We conclude from it that 
\begin{equation} \label{iniin}
  \rk\bar f\neq \rk f\;\;\Rightarrow\;\;\i_f\in [\B]:H_\B^\infty.
\end{equation}
We will use this statement later to justify the output $\{1\}$,
in case the condition in line 10 is satisfied. For now, assume
that $\rk\bar f=\rk f$.

Then, from Lemma~\ref{InitialBelongs}, we also have:
$$
H_{\bar f}\subset H_f\cdot H_\B^\infty+(\B).
$$
Since $f\in\C$, and due to invariants \ref{IB0} and \ref{IHB}, we thus obtain
$H_{\bar f}\subset H_\C^\infty+[\C]$. This means that 
\ref{IHB} is also preserved.

By definition of the algebraic remainder,
$$
f\in \left(\B\cup\left\{\bar f\right\}\right):H_\B^\infty.
$$
Note that line 6 results in adding $f$ to the set 
$\C\setminus\D$,
and this set is not changed elsewhere throughout the iteration.
Thus, at the end of the iteration \ref{IC} will hold.

Finally, as yet another consequence of Lemma~\ref{InitialBelongs},
$$
H_f\subset \left(H_{\bar f}+(\B)\right):H_\B^\infty.
$$
Taking into account that $\C\setminus\D$ does not change
other than in line 6, we thus obtain \ref{IHC} at the end of the 
iteration.

\item $f\in\delta\B^*_{\{m_i\}}$. As was shown
in Part 1, Case 2, $\B^0$ remains unchanged during the iteration
in this case. By \ref{IB0}, $f\in[\B^0]$, whence by 
definition of the algebraic remainder applied to line 8
we have
$$
\bar f\in \left(\B\cup\delta\B^0\cup f\right)\subset\left[\B^0\right].
$$
Thus, \ref{IB0} is preserved.

Next, according to \ref{IB} and \ref{IrkB}, all elements of $\B$,
and, hence, all elements of $\delta\B$, have distinct leaders.
In particular, if $\ld f\in\ld\delta\B^0$, then $f\in\delta\B^0$,
whence $\ld f\not\in\ld\left(\delta\B^0\setminus\{f\}\right)$.
In addition, since $f\in\delta\B$, and due to \ref{IB} and \ref{IrkB},
we have $\ld f\not\in\ld\B$. Altogether, 
$$
\ld f\not\in\ld\left(\B\cup\left(\delta\B^0\setminus\{f\}\right)\right).
$$
Thus, Lemma~\ref{InitialBelongs} (see also Remark~\ref{InitialRemark}) applies to the algebraic remainder computed in line 8, yielding (\ref{iniin}).
We will use this statement later to justify line 10, assuming
for now that $\rk\bar f=\rk f$.

From Lemma~\ref{InitialBelongs}, we also have:
$$H_{\bar f}\subset H_f\cdot H_{\B\cup\delta\B^0}^\infty+(\B\cup\delta\B^0)\subset H_f\cdot H_\B^\infty+[\B].$$
Since $f\in\delta\B$, and due to invariants \ref{IB0} and \ref{IHB}
and the fact that $H_f\subset H_\B$, we thus obtain
that \ref{IHB} is preserved at the end of the iteration.

Since set $\C\setminus\D$ remains unchanged during the iteration and 
set $\B$ is increased, invariants \ref{IC} and \ref{IHC} 
are automatically preserved. This concludes the study of Case 2.
\end{itemize}

Suppose now that $\rk \bar f<\rk f$ and apply the statement
(\ref{iniin}), which has been proved above in both cases.
According to \ref{IB0}, $\B\subset[\C]$, hence
$$[\B]:H_\B^\infty\subset [\C]:H_\B^\infty\subset [\C]:(H_\B\cup
H_\C)^\infty.$$
According to \ref{IHB}, $H_\B\subset H_C^\infty+[\C]$. Thus,
Lemma~\ref{l:Hubert} with $H=H_\C$, $K=H_\B$, and $I=[\C]$ yields
$[\C]:(H_\B\cup H_\C)^\infty=[\C]:H_\C^\infty$, whence
$[\B]:H_\B^\infty\subset [\C]:H_\C^\infty.$
In particular, keeping (\ref{iniin}) in mind, this implies that 
\begin{equation} \label{ifin1}
\i_f\in [\C]:H_\C^\infty.
\end{equation}
Due to \ref{IDC} for Case 1, or due to \ref{IHB} for Case 2, 
we also have that 
\begin{equation} \label{ifin2}
\i_f\in H_f\subset H_\C^\infty+[\C].
\end{equation}
Together (\ref{ifin1}) and (\ref{ifin2}) imply
$[\C]:H_\C^\infty=(1)$. This concludes the proof of correctness.
\end{enumerate}
\end{pf}

\subsection{Final algorithm and proof of the bound}\label{sec:finalalgorithm}

We are ready to present a modified version of the Rosenfeld-Gr\"obner
algorithm that satisfies the bound. The only place where the orders
of derivatives may grow is the pseudoreduction w.r.t. an autoreduced
set $\C$. Of course, only the orders of non-leading differential 
indeterminates may grow, while the orders of the leading ones decrease
as a result of reduction (or stay the same if the reduction turns out
to be algebraic, but then the orders of non-leading indeterminates
do not grow either). 

By associating different weights with leading
and non-leading indeterminates, we will achieve that the weighted
sum of their orders does not increase as a result of reduction. 
These weights come from the bound in the algorithm {\sf Differentiate\&Autoreduce}. 
If the set of leading indeterminates changes, so do the weights. 
However, if we estimate in advance the number of times the set of 
leading indeterminates can change throughout the algorithm, 
we can still obtain an overall bound on the orders.

For the original Rosenfeld-Gr\"obner algorithm, it is not that easy to 
carry out such an estimate, because some indeterminates may disappear 
and reappear again among the leading indeterminates of the
characteristic set $\C$. For example, 

\begin{exmp}
Let $F = \{y+z,$ $x,$ $x^2+z\},$ with the elimination ranking $x>y>z.$ 
\begin{itemize}
\item We choose its characteristic set as $\C := \{y+z,$ $x\}.$ 
\item The leading variables of $\C$ are $\{y,x\}.$ 
\item We put $\bar F := \Rem(F\setminus \C, \C) = \{z\}.$
\item $F_\mathrm{new} := \bar F \cup \C = \{z,$ $y+z,$ $x\}.$
\item As radical differential ideals:
$$
\left\{y+z,x,x^2+z\right\} = [z,y+z,x]:1^\infty\cap\left\{y+z,x,x^2+z,1\right\}.
$$
\item The new $\C =\{z,$ $x\}$ is computed from $F_\mathrm{new}$ and the leading variables have changed! 
\item \ldots
\bigskip
\item Finally,
$$
\left\{y+z,x,x^2+z\right\} = [z,y,x]:1^\infty = [z,y,x]
$$
and we see that the leaders $y$ and $x$ have come back.
\end{itemize}
\end{exmp}

\begin{exmp}
Let $F = \{zy,$ $x,$ $x^2+z\},$ with the elimination ranking $x>y>z.$ 
\begin{itemize}
\item We choose its characteristic set as $\C := \{zy,$ $x\}.$ 
\item The leading variables of $\C$ are $\{y,x\}.$ 
\item We put $\bar F := \Rem(F\setminus \C, \C) = \{z\}.$
\item $F_\mathrm{new} := \bar F \cup \C = \{z,$ $zy,$ $x\}.$
\item As radical differential ideals:
$$
\left\{zy,x,x^2+z\right\} = [z,zy,x]:z^\infty\cap\left\{zy,x,x^2+z,z\right\}.
$$
\item The new $\C =\{z,$ $x\}$ is computed from $F_\mathrm{new}$ and the leading variables have also changed! 
\item But the first component is trivial: $1 \in [z,zy,x]:z^\infty.$
\end{itemize}
\end{exmp}

The first situation can be remedied by properly relaxing the
requirement that $\C$ is autoreduced, while 
the second one can detected, after which further
computations in this branch of the Rosenfeld-Gr\"obner algorithm
are not necessary. As a result, we obtain an algorithm, in 
which, as long as an indeterminate appears among the leading 
indeterminates of the set $\C$, w.r.t. which we reduce,
it will stay there until the end.  

As mentioned above, we are going to replace the computation 
of the characteristic set by that of a weak d-triangular subset.
It is tempting to simply compute a weak d-triangular subset of 
the least rank, since this computation is inexpensive and it 
would give us the desired property that the leading indeterminates 
do not disappear. However, the termination of the algorithm 
is not guaranteed then. For example, take the system $F=\{x,xy\}$
in $\K\{x,y\}$, and let $x<y$. The weak d-triangular subset
of $F$ of the least rank is $F$ itself. Thus, we obtain a component
$\{x,xy\}:x^\infty=(1)$ and another component $\{x,xy,\i_{xy}\}$.
However, $\i_{xy}=x$, hence we arrive at the same set $F$ that
was given in the input, and the algorithm runs forever. 

The reason for the above behavior is that the initials of a 
weak d-triangular set $\C$, as opposed to an autoreduced set, 
need not be reduced w.r.t. $\C$. Thus by adding these 
initials we do not necessarily decrease the rank. The solution 
comes from the idea of \cite[Section 5]{Bou2}, \cite[Algorithm 6.11]{Dif}, and
\cite[Algorithm 4.1]{Imp} to construct the 
weak d-triangular set $\C$ gradually, so that each next polynomial $f$ 
to be added to $\C$ is reduced w.r.t. $\C$ (thus, we can also safely 
add the initial and separant of $f$ and guarantee that the rank 
decreases). In order to be able to construct the set $\C$
gradually, similarly to \cite{Dif}, we store it as a separate
component of the triples $(F,\C,H)\in U$.

The last modification that we are going to do is the replacement
of the differential pseudo-reduction w.r.t. $\C$ by 
the algebraic pseudo-reduction w.r.t. $\B$, which is computed
from $\C$ by Algorithm {\sf Differentiate\&Autoreduce}. 
As a result, we obtain Algorithm {\sf RGBound}. 

\begin{figure}
\begin{alg}{\sf RGBound}\label{RGBound}$(F_0,H_0)$\\
\begin{tabular}{l}
{\sc Input:} finite sets of differential polynomials
             $F_0\neq\varnothing$ and $H_0$, \\
\hphantom{\sc Input:} and a differential ranking\\
{\sc Output:} a finite set $T$ of regular systems such that \\
\hphantom{\sc Output:} $\qquad\qquad\{F_0\}:H_0^\infty=\bigcap\limits_{(\A,H)\in
  T}[\A]:H^\infty\;$ and \\
\hphantom{\sc Output:} $\qquad M(\A\cup H)\Le (n-1)!M(F_0\cup H_0)$ for $(\A,H)\in T.$\\
\\
\b $T:=\varnothing$, $\;\;\;U:=\{(F_0,\varnothing,H_0)\}$\\
\b {\bf while} $U\neq\varnothing$ {\bf do}\\
\b \b Take and remove any $(F,\C,H)\in U$\\
\b \b $f:=$ an element of $F$ of the least rank\\
\b \b $D:=\{C\in\C\;|\;\lv C=\lv f\}$\\
\b \b $G:=F\cup D\setminus\{f\}$\\
\b \b $\bar\C:=\C\setminus D\cup\{f\}$\\
\b \b $\B:=${\sf Differentiate\&Autoreduce}
      $\left(\bar\C, \left\{m_y(G\cup\bar\C\cup H)\;|\;y\in\lv\bar\C\right\}\right)$\\
\b \b {\bf if} $\B\neq\{1\}$ {\bf then}\\
\b \b\b $\bar F:=\algrem(G,\B)\setminus\{0\}$\\
\b \b\b $\bar H:=\algrem(H,\B)\cup H_\B$\\
\b \b\b {\bf if} $\bar F\cap\K=\varnothing$ {\bf and} $0\not\in\bar H$ {\bf then}\\
\b \b\b\b {\bf if} $\bar F=\varnothing$ {\bf then} 
            $T:=T\,\cup\,\left\{\left(\B^0,\bar H\right)\right\}$\\
\b \b\b\b\b {\bf else} $U:=U\cup\left\{\left(\bar F,\B^0,\bar H\right)\right\}$\\
\b \b\b\b {\bf end if}\\
\b \b\b {\bf end if}\\
\b \b {\bf end if}\\
\b \b {\bf if} $\s_f\not\in\K$ {\bf then}\\
\b \b \b $U:=U \cup\{(F\cup\{\s_f\},\C,H)\}$\\
\b \b \b {\bf if} $\i_f\not\in\K$ {\bf then} $U:=U\cup\{(F\cup\{\i_f\},\C,H)\}$
   {\bf end if}\\
\b \b {\bf end if}\\
\b {\bf end while}\\
\b {\bf return} $T$
\end{tabular}
\end{alg}
\end{figure}

In the proof of the bound, a key role is played by the quantity
$M_Z(F)$, which is defined for a finite set $F$ of differential
polynomials and a proper subset $Z\subsetneq Y$.
Assume that $|Z|=k<n$. As before, for a differential indeterminate $y\in Y$,
$m_y(F)$ denotes the highest order of a derivative of $y$ occurring
in $F$, or zero, if $y$ does not occur in $F$. Then
$$M_Z(F):=(n-k)\sum_{y\in Z} m_y(F)+\sum_{y\in Y\setminus Z} m_y(F).$$
We also recall the notation
$$M(F)=\sum_{y\in Y} m_y(F).$$

\begin{prop}\label{p:finalalgorithm} Algorithm~\ref{RGBound}  is correct and terminates
\end{prop}
\begin{pf}
We prove the following invariants of the {\bf while}-loop:
\begin{itemize}
  \item (I1) $\{F_0\}:H_0^\infty=\bigcap_{(F,\C,H)\in
    U}\{F\cup\C\}:H^\infty\cap\bigcap_{(\A,H)\in T}[\A]:H^\infty$
  \item For all $(F,\C,H)\in U$, 
  \begin{itemize}
     \item (I2) $\C$ is d-triangular,
     \item (I3) $F\neq\varnothing$ is reduced w.r.t. $\C$
     \item (I4) $H_\C\subset H$, 
     \item (I5) Let $l=|\lv\C|$. Then, if $l<n$,
          $$M_{\lv\C}(F\cup\C\cup H)\Le (n-1)\ldots(n-l)\cdot
          M(F_0\cup H_0),$$
          \hphantom{(I5)} otherwise $$M(F\cup\C\cup H)\Le (n-1)!\cdot M(F_0\cup H_0).$$
  \end{itemize}
\end{itemize}
The invariants hold for the initial triple $(F_0,\varnothing,H_0)$.
Assuming that they hold at the beginning of an iteration of the {\bf while} loop,
we will show that the invariants also take place at the end of the
iteration.

Let $(F,\C,H)$ be the triple taken and removed from $U$. 
Since $F\neq\varnothing$, we can compute an element $f\in F$ of the
least rank. Then $f$, as an element of $F$, is reduced w.r.t. $\C$. 
Applying \cite[Proposition 6.6]{Dif}, we have
$$\{F\cup\C\}:H^\infty=\{F\cup\C\}:(H\cup H_f)^\infty\cap
  \{F\cup\{\i_f\}\cup\C\}:H^\infty\cap
  \{F\cup\{\s_f\}\cup\C\}:H^\infty.$$
We note that, since $\rk\i_f<\rk f$ and $\rk\s_f<\rk f$,
polynomials $\i_f$ and $\s_f$ are, respectively, the elements of $F\cup\{\i_f\}$ and
$F\cup\{\s_f\}$ of the least rank (and, to repeat, their ranks are
less than the rank of the least element of $F$). Moreover, since in
the last two triples $(F\cup\{\i_f\},\C,H), (F\cup\{\s_f\},\C,H)$
only the first component has changed, invariants I2--I5 are preserved
for them. For the proof of invariant I1, it remains to show that 
\begin{equation} \label{toshow}
\{F\cup\C\}:(H\cup H_f)^\infty=\left\{\begin{array}{ll}
        \left[\B^0\right]:\bar H^\infty,& \bar F=\varnothing\\
        \left\{\bar F\cup\B^0\right\}:\bar H^\infty,&{\rm otherwise}.
        \end{array}\right.
\end{equation}
Given that $\C$ is d-triangular, the three assignments 
following the computation of $f$ ensure that $\bar\C$ is a weak 
d-triangular set of rank strictly less than $\C$,
because the polynomial $f$ is reduced w.r.t. $\C$ and we throw away (from $\C$)
all its elements with leading variables ``in conflict'' with the one of $f$.
We note that
$$G\cup\bar\C = (F\cup \D\setminus\{f\})\cup(\C\setminus \D)\cup\{f\} = F\cup\C.$$ 
Since $H_\C\subset H$, we also have
$H\cup H_f=H\cup H_{\bar\C}$. Therefore,
\begin{equation} \label{step1}
\{F\cup\C\}:(H\cup H_f)^\infty=\left\{G\cup\bar\C\right\}:\left(H\cup H_{\bar\C}\right)^\infty.
\end{equation}
Next, we use the properties of the set $\B$ ensured by Algorithm
{\sf Differentiate\&Autoreduce}. Since $H_\B\subset
H_{\bar\C}^\infty+\left[\bar\C\right]$, applying Lemma~\ref{l:Hubert} with $K=H_\B$, we obtain
\begin{equation} \label{step2}
\left\{G\cup\bar\C\right\}:\left(H\cup H_{\bar\C}\right)^\infty=
\left\{G\cup\bar\C\right\}:\left(H\cup H_{\bar\C}\cup H_{\B}\right)^\infty.
\end{equation}
The inclusions $\B\subset\left[\bar\C\right]$ and $\bar\C\subset[\B]:H_\B^\infty$
imply that 
\begin{equation} \label{step3}
\left\{G\cup\bar\C\right\}:\left(H\cup H_{\bar\C}\cup H_{\B}\right)^\infty=
\{G\cup\B\}:\left(H\cup H_{\bar\C}\cup H_{\B}\right)^\infty.
\end{equation}
Using the fact that  $H_{\bar\C}\subset
\left(H_\B^\infty+[\B]\right):H_\B^\infty$  (see Algorithm~\ref{DifferentiateAutoreduce})
and applying Lemma~\ref{l:Hubert} with $K=H_{\bar\C}$, we get
\begin{equation} \label{step4}
\{G\cup\B\}:\left(H\cup H_{\bar\C}\cup H_{\B}\right)^\infty=
\{G\cup\B\}:(H\cup H_{\B})^\infty.
\end{equation}
It follows from the definition of the algebraic pseudo-remainder ({\sf algrem}) that 
\begin{equation} \label{step5}
\{G\cup\B\}:(H\cup H_{\B})^\infty=
\{\bar F\cup\B\}:\bar H^\infty.
\end{equation}
Indeed, $\{G\cup\B\}:(H\cup H_\B)^\infty = \{\bar F\cup\B\}:(H\cup H_\B)^\infty.$
Take now any $f\in \left\{\bar F\cup\B\right\}:(H\cup H_\B)^\infty$. There exists 
$h \in (H\cup H_\B)^\infty$ such that $h\cdot f \in \left\{\bar F\cup\B\right\}.$
If $\bar h$ is a remainder of $h$ w.r.t. $\B$ then there exists $h' \in H_\B^\infty$
with $h'h - \bar h \in (\B).$ Hence, $$\bar h f \in \left\{\bar F\cup\B\right\}$$ and 
$$f \in \left\{\bar F\cup\B\right\}:\bar H^\infty.$$ The reverse inclusion is
done in a similar way.                                                                                       
Since $\B\subset\left[\B^0\right]$, we obtain that 
$\left\{\bar F\cup\B\right\}:\bar H^\infty = \left\{\bar F\cup\B^0\right\}:\bar H^\infty$.

The set $\B^0$ is d-triangular, its rank is equal to that of $\bar\C$, 
set $\bar H$ is partially reduced w.r.t. $\B^0$ and contains $H_\B^0$, 
and $\bar F$ is reduced w.r.t. $\B^0$. Moreover, if $\bar F=\varnothing$, we obtain 
the 
regular system $\left(\B^0,\bar H\right)$, 
which corresponds to 
the radical differential ideal 
$\left[\B^0\right]:\bar H^\infty=\left\{\bar F\cup\B^0\right\}:H_\B^\infty.$
Thus, we have proved \eqref{toshow} and also have demonstrated that
invariants I2--I4 hold for the triple $\left(\bar F,\B^0,\bar H\right)$. 

Termination of the algorithm is proved as follows. At each
iteration of the {\bf while}-loop, the triple $(F,\C,H)\in U$ is replaced
by at most three triples $\left(\bar F,\B^0,\bar H\right)$,
$(F\cup\{\i_f\},\C,H)$, and $(F\cup\{\s_f\},\C,H)$.

Define a relation $\prec$ on the set of all triples $(F,\C,H)$
satisfying I2--I4: let $(F',\C',H')\prec (F,\C,H)$ if and only if
either $\rk\C'<\rk\C$, or $\C'=\C$ and the element of the least
rank in $F'$ is strictly less than that in $F$. Then $\prec$ is
a lexicographic product of two well-orders, which is a well-order. 
We have shown that in the first triple we have $\rk\B^0<\rk\C$;
in the last two triples the second component $\C$ remains the same,
but the elements of the least rank of $F\cup\{\i_f\}$ and
$F\cup\{\s_f\}$ are strictly less than the element of $F$ of the least
rank. That is, each
of the three triples is less than $(F,\C,H)$ w.r.t.
the well-order $\prec$. 

Therefore, all triples computed by
the algorithm can be arranged in a ternary tree, in which 
$(F_0,\varnothing,H_0)$ is the root, and every path starting from the root
is finite. Let $\lambda$ be the maximal length of such a path. Then 
the number of vertices in the tree does not exceed $\sum_{i=0}^\lambda 3^i$.
Thus, the tree is finite, whence the algorithm terminates.

Finally, we show that invariant I5 holds for the triple 
$\left(\bar F,\B^0,\bar H\right)$. We assume that $|\lv\C|=l$.
Two cases are possible:
\begin{enumerate}
  \item $\lv f\in\lv\C$. Then $\lv\bar\C=\lv\C$ and for any finite set of
         polynomials $K$, if $l<n$, we have
   \begin{align} \label{trans0}
         M_{\lv\bar\C}(K)=M_{\lv\C}(K).
   \end{align}
  \item $\lv f\not\in\lv\C$. Then $\lv\bar\C=\lv\C\cup\{\lv f\}$ and
    $|\lv\bar\C|=l+1$. If $l+1<n$, we observe that 
    \begin{equation} \label{trans1}
        \begin{array}{l}
        M_{\lv\bar\C}(K)=\\
          \;\;  =(n-l-1)\sum\limits_{y\in\lv\bar\C}m_y(K)+\sum\limits_{y\not\in\lv\bar\C}m_y(K)=\\
          \;\;  =(n-l-1)\sum\limits_{y\in\lv\C}m_y(K)+(n-l-1)\cdot m_{\lv f}(K)+
             \sum\limits_{y\not\in\lv\bar\C}m_y(K)=\\
          \;\;  =(n-l-1)\sum\limits_{y\in\lv\C}m_y(K)+(n-l-2)\cdot m_{\lv f}(K)+
             \\
         \quad\quad +\left(m_{\lv f}(K) + \sum\limits_{y\not\in\lv\bar\C}m_y(K)\right)=\\
          \;\; =(n-l-1)\sum\limits_{y\in\lv\C}m_y(K)+\sum\limits_{y\not\in\lv\C}m_y(K)+
             (n-l-2)\cdot m_{\lv f}(K)\Le\\
          \;\; \Le(n-l)\sum\limits_{y\in\lv\C}m_y(K)+\sum\limits_{y\not\in\lv\C}m_y(K)+
             (n-l-2)\cdot M_{\lv\C}(K)=\\
          \;\; =(n-l-1)\cdot M_{\lv\C}(K)
        \end{array}
    \end{equation}
        (here we have used the fact that $m_{\lv f}(K)\Le
        M_{\lv\C}(K)$). 
\end{enumerate}
If $\lv\C<n$ and $|\lv\bar\C|=n$, we simply note that 
    \begin{equation} \label{trans2}
       M(K)\Le M_{\lv\C}(K).
    \end{equation}

Assume for simplicity that 
$$\ld\bar\C=\left\{y_1^{(d_1)},\ldots,y_k^{(d_k)}\right\},$$      
where $k=l$ or $k=l+1$. 
Since all derivatives of 
$y_i,$ $1 \Le i \Le k,$ presented in $F\cup\B\cup H$ of order 
greater than $d_i$ can be found among $\rk\B$, and since 
the elements of $\bar F$ and $\bar H\setminus H_\B$ are algebraic
pseudo-remainders of $G$ and $H$ w.r.t. $\B$, we have
\begin{equation} \label{ineqR}
m_i\left(\bar F\cup\B^0\cup\bar H\right) \Le \left\{
\begin{array}{ll}
d_i,& 1\Le i\Le k\\
m_i(G\cup\B\cup H),& k<i\Le n.
\end{array} \right.
\end{equation}
Also, recall that $\B$ satisfies the inequality 
(see \eqref{BInequality})
\begin{equation} \label{ineqB}
m_i(\B)\Le m_i\left(G\cup\bar\C\cup H\right)+\sum\limits_{j=1}^k (m_j\left(G\cup\bar\C\cup H\right)-d_j),\ k<i\Le n.
\end{equation}
Combining (\ref{ineqR}) and (\ref{ineqB}), we obtain that
\begin{align*}
M_{\lv\bar\C}\left(\bar F\cup\B^0\cup\bar H\right) &= (n-k)\sum_{i=1}^kd_i + \sum_{i=k+1}^n m_i\left(\bar F\cup\B^0\cup\bar H\right) \Le\\
&\Le (n-k)\sum_{i=1}^kd_i + \sum_{i=k+1}^n m_i(G\cup\B\cup H) \Le\\
&\Le (n-k)\sum_{i=1}^kd_i + \sum_{i=k+1}^n m_i\left(G\cup\bar\C\cup H\right) +\\
&\quad\ + (n-k)\sum_{j=1}^k(m_j\left(G\cup\bar\C\cup H\right)-d_j)=\\
&= (n-k)\sum_{i=1}^km_i\left(G\cup\bar\C\cup H\right) + \sum_{i=k+1}^n m_i\left(G\cup\bar\C\cup H\right) =\\
&= M_{\lv\bar\C}\left(G\cup\bar\C\cup H\right)
\end{align*}
and if $k = n$ then
$$
M\left(\bar F\cup\B^0\cup\bar H\right) = \sum_{i=1}^n d_i + 0 = M\left(G\cup\bar\C\cup H\right)
$$
because $\rk\bar\C = \rk\B^0.$
Thus, 
\begin{equation} \label{ineqtrans1}
\begin{array}{ll}
M_{\lv\bar\C}\left(\bar F\cup\B^0\cup\bar H\right)\Le M_{\lv\bar\C}\left(G\cup\bar\C\cup H\right),&k<n\\
M\left(\bar F\cup\B^0\cup\bar H\right)\Le M\left(G\cup\bar\C\cup H\right),&k=n.
\end{array}
\end{equation}
Now, applying \eqref{trans0}, \eqref{trans1}, or \eqref{trans2} with
$K=G\cup\bar\C\cup H=F\cup\C\cup H$, we get

\begin{equation}\label{ineqtrans2}
\begin{array}{rlr}
M_{\lv\bar\C}\left(\bar F\cup\B^0\cup\bar H\right)&\Le M_{\lv\C}(F\cup\C\cup H),&l=k<n\\
M_{\lv\bar\C}\left(\bar F\cup\B^0\cup\bar H\right)&\Le (n-l-1)M_{\lv\C}\left(\bar F\cup\B^0\cup\bar H\right)\Le\\ 
& \Le (n-l-1)\cdot M_{\lv\C}(F\cup\C\cup H),& l<k<n\\
M\left(\bar F\cup\B^0\cup\bar H\right)&\Le M\left(G\cup\bar\C\cup H\right)=\\ &= M(F\cup\C\cup H)\Le\\
& \Le M_{\lv\C}(F\cup\C\cup H),& l<k=n\\
M\left(\bar F\cup\B^0\cup\bar H\right)&\Le M\left(G\cup\bar\C\cup H\right)=\\ &= M(F\cup\C\cup H),&l=k=n.
\end{array}
\end{equation}

By taking into account the fact that invariant I5 holds for the triple $(F,\C,H)$,
we thus obtain this invariant for the triple $\left(\bar F,\B^0,\bar H\right)$. 

To conclude the proof of the bound for the output regular systems
$\left(\B^0,\bar H\right)$, we note that it is already given by the invariant I5 when $k=n$,
while in case $k<n$ we 
use inequality \eqref{trans2}:
$$
M\left(\B^0\cup\bar H\right) \Le M_{\lv\B^0}\left(\B^0\cup\bar H\right)\Le (n-1)!\cdot M(F_0\cup H_0).
$$
\end{pf}

\subsection{Reduction-independent algorithm}
In Algorithm~\ref{RGBound} we had to be very careful in the reduction
process. The idea was to emulate differential reductions by doing enough
differentiations first and then applying purely algebraic reduction.
We take care of the orders of derivatives in the first process and do
not need to worry about them during the second purely algebraic step. 
Let us find out why such two-step procedure was necessary. If we reduce 
w.r.t. an arbitrary d-triangular set, the result of reduction depends on 
the choice of the reduction path.
\begin{exmp} Consider the following differential chain
$$
\C = x(x-1),\ (x-1)y,\ z+y+tx
$$ with the elimination ranking $t<x<y<z$
and the differential polynomial $$f = z'+y'.$$ We can reduce $f$ w.r.t.
$\C$ in many different ways and the remainders are very different:
\begin{enumerate}
\item \begin{equation*}
\begin{CD}
z'+y' @>{z'+y'+t'x+tx'}>> t'x+tx' @>{(2x-1)x'}>>t'(2x-1)x = 2t'x^2-t'x@>{}>> \\
@>{x^2-x}>>t'x =: f_1\\
\end{CD}
\end{equation*}
\item
\begin{equation*}
\begin{CD}
z'+y' @>{(x-1)y' + x'y}>> (x-1)z' - x'y @>{(x-1)y}>> (x-1)^2z'@>{}>> \\
@>{z'+y'+t'x+tx'}>>(x-1)^2(y'+t'x+tx')@>{}>> 0 =: f_2.\\
\end{CD}
\end{equation*}
\end{enumerate}
\end{exmp}

We see that the remainder $f_1$ depends on the variable $t'$ that
is not in both $f_2$ and $\C.$
So, the reason for these so different answers is that the set $\C$
has a non-invertible initial. Speaking informally, if $\C$ is partially autoreduced and 
its initials and separants are invertible, then the result of reduction is more or less
uniquely determined. More precisely, one can show that all 
results of reduction of a polynomial w.r.t. a d-triangular set with 
invertible initials and separants lie in a fixed N\"otherian ring of algebraic
polynomials. In particular, if one of the results of reduction satisfies a certain bound
on the order of its derivatives, then any other result of reduction will satisfy this
bound as well. 

Since we are not in position of reducing w.r.t. a set with invertible initials
and separants, we are going to state precisely and prove a slightly weaker statement.
Within the scope of this section, let us 
 call polynomial $g$ a {\it differential 
remainder} of polynomial $f$ w.r.t. $\C$, if $g$ is reduced w.r.t. 
$\C$ and there exists $h\in H_\C^\infty$ such that 
$$hf-g\in[\C]:H_\C^\infty.$$

\begin{prop}\label{p:SmallRing} Let $\C$ be a coherent\footnote{The
adjective ``coherent'' makes the statement valid
in presence of partial derivatives; in the ordinary case, it can be
ignored.}
d-triangular set of differential polynomials, $f$ a differential polynomial, 
and $g$ a differential remainder of $f$ w.r.t. $\C$.
Let $X$ be the set of derivatives present in $\C$ and $g$.
Let $\bar g$ be another differential remainder of $f$ w.r.t. $\C$.
Assume that $\bar g$ is not in $\K[X]$, i.e., it admits a representation 
$\bar g = a_ku^k + \ldots + a_0$, where $u\not\in X$
and $a_0,\ldots,a_k$ are free of $u$. Then $$\bar g - a_0 \in (\C):H_\C^\infty.$$
In particular, $a_0$ is also a differential remainder of $f$ w.r.t. $\C$.
\end{prop}
\begin{pf} 
Since $g$ and $\bar g$ are differential remainders of $f$ w.r.t. $g$, they 
are both reduced w.r.t. $\C$, and there
exist $h,\bar h\in H_\C^\infty$ such that 
$$h f - g \in [\C]:H_\C^\infty,\quad\bar h f - \bar g \in [\C]:H_\C^\infty.$$

Consider the differential polynomial
$$
\bar f := \bar h(hf-g) - h(\bar h f -\bar g) = h\bar g  - \bar h g\in [\C]:H_\C^\infty.
$$
Since $\C$ is a coherent d-triangular set, ideal $[\C]:H_\C^\infty$
is regular. The polynomial $\bar f$ is partially reduced w.r.t. $\C.$ 
Therefore, by the Rosenfeld Lemma $\bar f\in(\C):H_\C^\infty$. 

We have $$\bar f =  (h\cdot a_k)u^k +\ldots + (h\cdot a_0 - \bar h\cdot g)$$
with  $\bar h\cdot g$ contributing only to $a_0,$ because it does not
depend on $u.$ Since $u$ does not appear in $\C,$ the fact
that $\bar f\in(\C):H_\C^\infty$ implies 
that every coefficient of $\bar f$ belongs to this ideal. In
particular, $h\cdot a_k$ belongs
to $(\C):H_\C^\infty,$  whence $a_i \in (\C):H_\C^\infty,$ $1 \Le i \Le k.$ Thus, 
$$\bar g - a_0=a_ku^k+\ldots+a_1u \in (\C):H_\C^\infty.$$
\end{pf}

We are going to apply the above Proposition as follows. Let $\C$ and $f$
be as in its statement. 
Suppose we know that there exists a differential remainder $g$ of 
$f$ w.r.t. $\C$ that satisfies a certain bound $b$ on the order of derivatives
occurring in it. We emphasize that we do not need to know $g$, 
the fact of its existence is sufficient. Compute {\it any} differential remainder
$\bar g$ of $f$ w.r.t. $\C$. Then, if $\bar g$ does not satisfy the bound $b$,
it must contain a derivative $u$ that does not satisfy this bound. 
By Proposition~\ref{p:SmallRing}, the free coefficient of $\bar g$, when 
viewed as a polynomial in $u$, is also a differential remainder of $f$ w.r.t. 
$\C$. Replace $\bar g$ by its free coefficient; continue such replacements
until $\bar g$ satisfies the bound $b$. This yields an efficient procedure
that computes a differential remainder satisfying the bound:

\begin{figure}[!tbh]
\begin{alg}{\sf
    Truncate}\label{Tailpart} $(f, \{p_i\})$\\
\begin{tabular}{l}
{\sc Input:} 
a differential polynomial $f$ and numbers $p_i\Ge 0$\\
{\sc Output:} {\bf truncation} of $f$, i.e., the sum of those terms of $f$ that \\
\hphantom{\sc Output:} belong 
to the polynomial ring $R=\K[\ldots,y_i,\dots,y_i^{(p_i)},y_{i+1},\ldots]$\\
\b Let $f=\alpha_1+\ldots+\alpha_q$, where $\alpha_i$ are differential monomials\\
\b $g:=0$\\
\b {\bf for} $i:=1$ {\bf to} $q$ {\bf do}\\
\b \b {\bf if} $\alpha_i\in R$ {\bf then} $g:=g+\alpha_i$\\
\b {\bf end for}\\
\b {\bf return} $g$
\end{tabular}
\end{alg}
\end{figure}

We have proved the following
\begin{thm} \label{Truncation}
   Let $\C$ be a coherent d-triangular set of differential polynomials, 
   and let $f$ be a differential polynomial. Let $p_i\Ge m_i(\C)$, $i=1,\ldots,n$.
   Assume that there exists
   a differential remainder of $f$ w.r.t. $\C$, which contains no derivatives
   of differential indeterminate $y_i$ of order greater than $p_i$, $i=1,\ldots,n$.
   Let $g$ be any differential remainder of $f$ w.r.t. $\C$. Then 
   {\sf Truncate}$(g,\{p_i\})$ is a differential remainder of $f$ w.r.t. $\C$,
   in which the order of every differential indeterminate $y_i$ does not exceed
   $p_i$. 
\end{thm}

We are going to modify Algorithm~\ref{RGBound}, so that there is no necessity
to perform differential pseudo-reduction in two steps, via prolongation and
purely algebraic reduction. In the new Algorithm~\ref{RGBoundRI}, it is assumed that 
procedure $\Rem$ computes any differential remainder in the above sense. 
The key idea is the following: whenever we find a differential remainder w.r.t. $\D$ 
that does not satisfy the expected bound $b$ (computed by Algorithm~\ref{RGBoundRI}), 
by Theorem~\ref{Truncation} we can simply truncate this remainder. In order 
to be able to apply Theorem~\ref{Truncation}, we are going to prove
the existence of a differential remainder satisfying $b$. In fact, we 
know that sets $\B$, $\bar F$, and $\bar H$ computed in Algorithm~\ref{RGBound}
satisfy $b$; it remains to be shown that one can obtain differential remainders w.r.t. $\D$,
given the elements of $\B$, $\bar F$, and $\bar H$.
Note that we may assume
$\rk\D=\rk\bar\C$ (at the end of the {\bf for}-loop), since otherwise all results of truncations 
are discarded by Algorithm~\ref{RGBoundRI}. 

To justify truncations in the {\bf for}-loop of Algorithm~\ref{RGBoundRI}, 
we consider
$$
\B = {\sf Differentiate\&Autoreduce}\left(\bar\C, \left\{m_y\left(\bar G\right)\:|\: y \in \lv\bar\C\right\}\right).
$$
and show that, at the beginning of each iteration, there exist $B\in\B$ and $h\in H_\D^\infty$
such that $hB$ is a differential remainder of $C$ w.r.t. $\D$.
This statement is a consequence of the following expanded invariant of the {\bf for}-loop, 
which we are going to prove by induction 
on the number of iterations. Let 
$$\B_{<C}=\{f\in\B\;|\;\ld f<\ld C\},$$
$B = \algrem(C,\B_{<C}),$ $E = \Rem(C,\D),$ and $D =$ {\sf Truncate} $(E,b).$
Then
\begin{align*}
h'\cdot C - hB \in [\D]:H_\D^\infty,\\
B\in[\D\cup\{D\}]:H_\D^\infty,\\
H_B \subset (H_D^\infty + [\D]):H_\D^\infty,\\
\end{align*}
for some $h,h' \in H_\D^\infty.$ 

The inductive base holds, since at the end of the first iteration we have
$B=E=D=C$ and $\D=\{C\}$. 

For the inductive step, we have: $h_1\cdot C - B \in \left(\B_{<C}\right)$ 
for some $h_1 \in H_{\B_{<C}}^\infty.$
By the inductive assumption $[\B_{<C}] \subset [\D]:H_\D^\infty.$
Hence, 
$$
h_1\cdot C - B \in [\D]:H_\D^\infty.
$$
Also by the inductive assumption, $h_1 \in (H_\D^\infty + [\D]):H_\D^\infty.$ 
This means that there exist $h \in H_\D^\infty,$ $h' \in H_\D^\infty$ such
that $h\cdot h_1 - h' \in [\D].$ Thus, $$h'\cdot C - h\cdot B \in [\D]:H_\D^\infty.$$

By definition of (algebraic) pseudo-remainder, we have
$$
B\in\left(\B_{<C}\cup\{C\}\right),\quad C\in (E+[\D]):H_\D^\infty.
$$
By Lemma~\ref{InitialBelongs}, taking into account the assumption $\rk\D=\rk\bar\C$, we have:
$$
H_B\subset H_C\cdot H_{\B_{<C}}^\infty+\left(\B_{<C}\right),\quad H_C\subset (H_E+[\D]):H_\D^\infty.
$$
By Proposition~\ref{p:SmallRing}, $E\in D+(\D):H_\D^\infty$. By modifying slightly 
the proof of Lemma~\ref{InitialBelongs}, we will show that this implies
$H_E\subset H_D+(\D):H_\D^\infty$. Indeed, using the assumption $\rk\D=\rk\bar\C$ (which holds at
the end of the {\bf for}-loop), we
obtain $\rk D=\rk C=\rk E$; since all leading differential indeterminates in $\bar\C$ are distinct, this, 
in particular, implies that 
$$
v=\ld D=\ld E\not\in\ld\D.
$$ 
Now let $f_1,\ldots,f_k$ be any generators of the ideal $(\D):H_\D^\infty$, so that we have
$$
E-D=\sum_{i=1}^k \alpha_i f_i.
$$
By viewing the above equality as one between two polynomials in $v$ and noting
that $f_i$ do not involve $v$, we immediately
obtain that $\i_E-\i_D\in (\D):H_\D^\infty$ and $\s_E-\s_D\in (\D):H_\D^\infty$.

Combining the above statements, we obtain the required invariants at the 
end of the iteration:
$$
B\in\left(\B_{<C}\cup\{C\}\right)\subset ([\D]:H_\D^\infty\cup(E+[\D]):H_\D^\infty)\subset[\D\cup\{D\}]:H_\D^\infty
$$
and 
$$
H_B\subset H_C\cdot H_{\B_{<\C}}^\infty+\left(\B_{<C}\right)\subset (H_E+[\D]):H_\D^\infty+[\D]:H_\D^\infty
\subset
(H_D+[\D]):H_\D^\infty.
$$

The truncations applied in Algorithm~\ref{RGBoundRI} to compute sets $\bar F$ and $\bar H$ 
are justified by showing that differential remainders of $G$ and $H\cup H_f$ w.r.t. $\D$ 
that satisfy the bound $b$ exist and can be similarly obtained from the elements of sets 
$\bar F$ and $\bar H$ computed by Algorithm~\ref{RGBound}. We omit these details.

\begin{figure}
\begin{alg}{\sf RGBound-Reduction-Independent}\label{RGBoundRI}$(F_0,H_0)$\\
\begin{tabular}{l}
{\sc Input:} finite sets of differential polynomials
             $F_0\neq\varnothing$ and $H_0$, \\
\hphantom{\sc Input:} and a differential ranking\\
{\sc Output:} a finite set $T$ of regular systems such that \\
\hphantom{\sc Output:} $\qquad\qquad\{F_0\}:H_0^\infty=\bigcap\limits_{(\A,H)\in
  T}[\A]:H^\infty\;$ and \\
\hphantom{\sc Output:} $\qquad M(\A\cup H)\Le (n-1)!M(F_0\cup H_0)$ for $(\A,H)\in T.$\\
\b $T:=\varnothing$, $\;\;\;U:=\{(F_0,\varnothing,H_0)\}$\\
\b {\bf while} $U\neq\varnothing$ {\bf do}\\
\b \b Take and remove any $(F,\C,H)\in U$\\
\b \b $f$ an element of $F$ of the least rank\\
\b \b $D:=\{C\in\C\;|\;\lv C = \lv f\}$\\
\b \b $G:=F\cup D\setminus \{f\}$\\
\b \b $\bar\C:=\C\setminus D\cup\{f\}$\\
\b \b $\bar G:=G\cup\bar\C\cup H$\\
\b \b $b := \left\{m_y\left(\bar G\right)\:\big|\:y\in \lv\bar C\right\}\cup
\left\{m_z\left(\bar G\right)+\sum\limits_{y\in\lv\bar C}\left(m_y\left(\bar G\right)-m_y\left(\ld\bar\C\right)\right)\:\big|\: z \notin \lv\bar\C\right\}$\\
\b \b $\D := \varnothing$\\
\b\b {\bf for} $C \in \bar\C$ increasingly {\bf do}\\
\b \b \b  $\D := \D \cup \left\{{\sf Truncate} \left(\Rem\left(C,\D\right),b\right) \right\}$ \\
\b \b {\bf if} $\rk\D = \rk\bar\C$ {\bf then}\\
\b\b \b $\bar F:=$ {\sf Truncate} $\left(\drem(G,\D)\setminus\{0\},b\right)$\\
\b\b \b $\bar H:=$ {\sf Truncate} $\left(\drem(H\cup H_f,\D)\cup H_{\D},b\right)$\\
\b\b \b {\bf if} $\bar F\cap\K=\varnothing$ {\bf and} $0\not\in\bar H$ \\
\b\b\b\b {\bf then} $U :=  U\cup\left\{\bar F, \D, \bar H \right\}$ \\
\b\b\b\b {\bf else} $T:=T\,\cup\,\{(\D,\bar H)\}$ \\
\b\b\b {\bf end if}\\
\b\b{\bf end if}\\
\b \b {\bf if} $\s_f\not\in\K$ {\bf then}\\
\b \b \b $U:=U \cup\{(F\cup\{\s_f\},\C,H)\}$\\
\b \b \b {\bf if} $\i_f\not\in\K$ {\bf then} $U:=U\cup\{(F\cup\{\i_f\},\C,H)\}$
   {\bf end if}\\
\b \b {\bf end if}\\
\b {\bf end while}\\
\b {\bf return} $T$
\end{tabular}
\end{alg}
\end{figure}

\begin{prop} Algorithm~\ref{RGBoundRI} is correct and satisfies the bound.
\end{prop}
\begin{pf}
The proof of correctness, termination, and bound for 
Algorithm~\ref{RGBoundRI} is based on the same invariants
of the {\bf while}-loop that were used for Algorithm~\ref{RGBound}. The
only new step we make is the {\sf Truncate} algorithm whose application
 is justified above. 
\end{pf}

\def\T{\mathbb T}
\def\F{\mathbb F}

\section{Transformation of characteristic sets of prime differential
  ideals} \label{sec:prime}

As above, let $\K\{Y\}$ be a ring of ordinary differential polynomials
in $n$ indeterminates with the derivation $\delta$. 
Let $\C$ be a characteristic set of a prime differential ideal $I$
in $\K\{Y\}$ w.r.t. a ranking $\le$. We propose an algorithm 
that computes a characteristic set of $I$ w.r.t. any other ranking
$\le'$ {\it algebraically}. More precisely, using a bound on 
the orders of derivatives occurring in the canonical characteristic 
set $\D$ of $I$ w.r.t. the target ranking, we find a sufficient {\it
  differential prolongation} of $\C$ (described below), which 
defines a prime algebraic sub-ideal $\bar I$ in $I$ containing 
$\D$. After that, it remains to compute an algebraic 
characteristic set of $\bar I$ w.r.t. the target ranking and 
extract from it a differential characteristic set of $I$. 

\subsection{A bound for characteristic sets of prime differential
ideals} \label{sec:primebound}

First, given a characteristic set $\C$ of a prime differential 
ideal $I$ w.r.t. an arbitrary ranking $\le$, we would like
to obtain a bound on the orders of derivatives occurring in 
a characteristic set of $I$ w.r.t. another given ranking $\le'$.
For $\le$ orderly and $\le'$ arbitrary, such a bound is given 
in \citep{canonical}. If $\le$ is not orderly, we first obtain
a bound for the orders of the elements of an orderly characteristic
set $\D$ of $I$, and then apply the bound from \citep{canonical}.

Indeed, $\D$ can be computed from $\C$ with the help
of the Rosenfeld-Gr\"obner algorithm applied to the system
$F_0=\C$, $H_0=H_\C$ (where the initials and separants of $\C$
in $H_\C$ are taken w.r.t. $\le$). Since $I$ is prime,
one of the regular components $(A,H)$ computed by the Rosenfeld-Gr\"obner
algorithm will coincide with $I$, and the characteristic set 
of the corresponding regular ideal $[A]:H^\infty$ w.r.t. $\le'$
can be extracted from the lexicographic Gr\"obner basis of the
algebraic ideal $(A):H^\infty$ via the algorithm
given in \cite[Theorem 6]{Bou1}. A more
efficient algorithm, which uses the fact that the given ideal
is prime and thus avoids the computation of redundant regular
components, is presented in \cite{PARDI}. 

Let $M$ be the maximal order of derivatives occurring in $\C$.
The only place where the Rosenfeld-Gr\"obner algorithm differentiates 
polynomials is the computation of differential pseudo-remainders.
However, for an orderly ranking, the order of a polynomial cannot
increase as a result of pseudo-reduction. Thus, the orders of derivatives
occurring in the characteristic set $\D$ do not exceed $M$.
In fact, the same applies to any other characteristic set of $I$ 
w.r.t. the same orderly ranking: the leading derivatives of all characteristic 
sets of $I$ w.r.t. the same ranking coincide, and the orders of non-leading 
derivatives occurring in a polynomial $f$ cannot exceed the order 
of the leader of $f$ w.r.t. an orderly ranking. 

Now we will use the following 
\begin{lem} \label{l:cardcs}
  The number of elements in a characteristic set $\C$ of a prime
  differential ideal $I$ in the ring of ordinary differential
  polynomials  $\K\{y_1,\ldots, y_n\}$ does not depend on the ranking.
\end{lem}
\begin{pf} If $d$ is a differential dimension of $P$ then the number of elements of $\C$ is equal to $n-d$ by \cite[Theorem 4.11]{Resolvent} which
does not depend on a choice of a differential ranking.
\end{pf}

\begin{rem}
  The above lemma does not hold in the partial differential case.
  For example (borrowed from \cite{PARDI}), a characteristic set of the prime differential ideal 
  $$[u_x^2-4u,\;u_{xy}v_y-u+1,\;v_{xx}-u_x]$$ in $\K\{Y\}$ with 
  derivations $\Delta=\{\partial/\partial x,\partial/\partial y\}$  
  may have 3 or 4 elements, depending on the ranking. 
\end{rem}

For the above example, It takes a while to compute the characteristic 
set of the ideal w.r.t. the elimination ranking $u > v$ using the 
Rosenfeld-Gr\"obner algorithm in Maple (see \cite{ucs}).
Consider another example that requires less computational efforts.

\begin{exmp} Consider the following prime differential ideal: $$P = [u_{yy},\;v_{xx}+y\cdot u_x+u].$$ This set of generators forms
a characteristic set of $P$ w.r.t. the elimination ranking with $v > u.$
However, if we change the ranking to $u > v$, then the following set 
containing $3$ elements will be a characteristic set of $P:$
\begin{align*}
&v_{xxyyy},\\ 
&y^2\cdot v_{xxxxyy}- 2y\cdot v_{xxxxy} + 2y\cdot v_{xxxyy}  + 2v_{xxxx} - 2v_{xxxy} + v_{xxyy},\\
&2u - y^3\cdot v_{xxxyy}+2y^2\cdot v_{xxxy} -2y\cdot v_{xxx} +2v_{xx}.
\end{align*}
\end{exmp}

Applying Lemma~\ref{l:cardcs}, we obtain the following bound on the
{\it order} of $I$ (see \cite{canonical}):
\begin{equation} \label{bound}
 \ord I:=\sum_{D\in\D}\ord D\Le |\C|\cdot\max_{C\in\C}\ord C.
\end{equation}
This bound is likely to be non-optimal. It is possible that the results
of \cite[Chapter VII]{Rit}, together with Lemma~\ref{l:cardcs}, imply
the following bound, which is better: 
let $m_1\Ge m_2\Ge\ldots \Ge m_n$ be the numbers $m_y(\C)$, $y\in Y$,
arranged in non-increasing order, then $$\ord I\Le \sum_{i=1}^{|\C|}
m_i.$$ For this bound, which so far is a conjecture, 
one needs to verify that Ritt's proof holds for
non-elimination rankings and also adapt it for 
ideals specified by characteristic sets, rather
than sets of generators. 

According to \citep{canonical}, the orders of derivatives occurring
in the canonical characteristic set of $I$ w.r.t. any ranking
do not exceed the order of $I$. Thus, the number
$$M_1=|\C|\cdot\max_{C\in\C}\ord C$$
bounds the orders of derivatives occurring in the canonical characteristic
set of $I$ w.r.t. any (not necessarily orderly) target ranking $\le'$.

We note that the bound $(n-1)!\cdot M(\C)$ obtained in Section~\ref{sec:finalalgorithm}
is also a bound for the orders of derivatives occurring in the 
characteristic set of $I$ w.r.t. $\le'$ computed by the 
Rosenfeld-Gr\"obner algorithm. In fact, invariant I5 in the proof of 
Proposition~\ref{p:finalalgorithm}, together with Lemma~\ref{l:cardcs}, yields a better bound 
$$M_2=\frac{(n-1)!}{(n-|\C|-1)!}\cdot M(\C).$$
In most cases, $M_2>M_1$, but in some, especially for small values of $n$, it may
happen that $M_2<M_1$. This again suggests that none of the two bounds is
optimal. Leaving the important problem of obtaining an optimal bound for future
research, we summarize the bounds obtained so far in the following
\begin{lem} \label{l:boundprime}
  Let $\C$ be a characteristic set of an ordinary prime differential 
  ideal $I$ w.r.t. a ranking $\le$. Then $\ord I$ and 
  the orders of derivatives occurring 
  in the canonical characteristic set of $I$ w.r.t. another ranking
  $\le'$ do not exceed 
$$
M_\C:=\min(M_1,M_2)=\min\left(|\C|\cdot\max_{C\in\C}\ord C,\frac{(n-1)!}{(n-|\C|-1)!}\cdot M(\C)\right).
$$
\end{lem}

\subsection{Differential prolongation: the prime case}

Assume that $\ld_\le\C=\left\{y_1^{(d_1)},\ldots,y_k^{(d_k)}\right\}$.
Let $m_i=M_\C$, $1\Le i\Le k$.
Compute the set
$$\A={\sf Differentiate\&Autoreduce}\left(\C,\{m_i\}_{i=1}^k\right)$$
(for the algorithm {\sf Differentiate\&Autoreduce}, see 
Section~\ref{sec:da} above). 

Let $\D$ be the canonical characteristic set of $I$ w.r.t. $\le'$. 
Every polynomial in $\D$, as an element of $I$, reduces 
w.r.t. $\C$ and $\le$ to zero. Since the orders of 
derivatives occurring in $\D$ do not exceed $M_\C$,
every polynomial in $\D$ algebraically reduces to zero w.r.t. $\A$. 
That is, $\D\subset(\A):H_\A^\infty$.

The algebraic ideal $\bar I=(\A):H_\A^\infty$ is equal to 
the intersection of $I$ with the ring
$R=\K\left[\Theta Y\setminus\Theta\ld_\le\C\cup\ld_\le\A\right].$
Indeed, $\A\subset R$. Vice versa, every element of $I\cap R$
algebraically reduces w.r.t. $\A$ to zero and therefore
belongs to $(\A):H_\A^\infty$. 

Since $I$ is prime, so is $\bar I$. 
Applying one of the existing efficient algorithms (for instance, see \cite{PARDI} or \cite{MarcMaple}) to 
the set $\A$, we compute the canonical algebraic characteristic set $\B$
of $\bar I$ w.r.t. the target ranking $\le'$. We know that
the algebraic ideal $\bar I$ contains the canonical characteristic
set $\D$ of the differential ideal $I$ w.r.t. $\le'$.
In the following section, we will show that, in fact, $\D\subseteq\B$.

\subsection{Extracting a differential characteristic set}

The following two lemmas hold in the partial differential case. 
We assume that a ranking is fixed. 

\begin{lem} \label{l:rkweakdtr}
  Let $\K\{Y\}$ be a ring of partial differential polynomials, and let $K$ be 
  an arbitrary subset of $\K\{Y\}\setminus\K$.

  Let $\C$ be a differential characteristic set of $K$ and 
  $\A$ an algebraic characteristic set of $K$. 
  Let $\T$ be a weak d-triangular subset of $\A$ of the least rank.
  Then $\rk\T\le\rk\C$. 
\end{lem}
\begin{pf}
  Suppose that a polynomial $f\in\C$ is differentially
  reduced w.r.t. $\T$. Then, since $\T$ is a weak d-triangular subset 
  of $\A$ of the least rank, $f$ is algebraically reduced w.r.t. $\A$.
  Due to the fact that $\A$ is an algebraic characteristic set of $K$,
  we have $f=0$, contradiction. Thus, no element of $\C$ is
  differentially reduced w.r.t. $\T$, which implies that $\rk\T\le\rk\C$.
\end{pf}

\begin{lem} \label{l:algdifcs}
  Let $I$ be a prime differential ideal, let $\C$ be
  the canonical characteristic set of $I$, and let 
  $J=I\cap\K[V]$, where $V\subset\Theta Y$,
  be an algebraic ideal containing $\C$. Then 
  the canonical algebraic characteristic set $\D$ of $J$ 
  contains $\C$; more precisely, $\C$ is the weak d-triangular
  subset of $\D$ of the least rank. 
\end{lem}  
\begin{pf}
  Since $\D$ is triangular, its weak d-triangular subset
  of the least rank is unique. Let $\T$ be the weak d-triangular
  subset of $\D$ of the least rank. 

  Since $\D$ is an algebraic characteristic set of the prime ideal $J$, we have
  $H_\D\cap J=\varnothing$. Moreover, $H_\D\subset\K[V]$, 
  therefore $H_\D\cap I=\varnothing$ and, hence, 
  $H_\T\cap I=\varnothing$. Since $\T\subset I$ and $I$ is prime, this implies 
  \begin{equation} \label{tincl}
     [\T]:H_\T^\infty\subset I.
  \end{equation}

  Let 
  $$\A=\{\drem(f,\T\setminus\{f\})\;|\;f\in\T\}.$$
  We have $\A\subset[\T]\subset I$; we will show that set $\A$ is 
  differentially autoreduced and $\rk\A=\rk\T$. 

  First, show that $\rk\A=\rk\T$. Indeed, suppose that for some
  $f\in\T$ and $g=\drem(f,\T\setminus\{f\})$, we have $\rk g<\rk f$. 
  Since $\T$ is weak d-triangular, $\ld f\not\in\Theta\ld(\T\setminus\{f\})$. 
  Thus, Lemma~\ref{InitialBelongs} applies and tells us that 
  $\i_f\in[\T]:H_\T^\infty$. Hence, according to \eqref{tincl},
  $\i_f\in I$. This contradicts with the fact that 
  $H_\T\cap I=\varnothing$. 

  Now, since $g$ is reduced w.r.t. $\T\setminus\{f\}$, 
  $\rk g=\rk f$, and $\rk\A=\rk\T$, 
  $g$ is also reduced w.r.t. $\A\setminus\{g\}$. That is,
  set $\A$ is autoreduced. 

  By Lemma~\ref{l:rkweakdtr}, $\rk\T\le\rk\C$. Therefore,
  $\rk\A\le\rk\C$. 
  Since $\A$ is an autoreduced subset of $I$, while $\C$
  is an autoreduced subset of $I$ of the least rank,
  we have $\rk\A\ge\rk\C$. Thus, $\rk\A=\rk\T=\rk\C$.

  Let $\bar\D=(\D\setminus\T)\cup\C$. Set $\bar\D$ is
  algebraically autoreduced, has the same rank as $\D$, and satisfies
  the requirements of canonicity: for every $f\in\bar\D$, 
  the initial of $f$ does not depend on the leaders of $\bar\D$,
  $f$ is monic and has no factors in $\K[N(\bar\D)]$, where 
  $N(\bar\D)=N(\D)=V\setminus\ld\D$ is the 
  set of non-leaders of $\D$ (or $\bar\D$). Since the canonical
  characteristic set is unique, we have $\bar\D=\D$ and 
  $\C=\T$. This concludes the proof. 
\end{pf}

Returning to the notation from the previous section and applying the above
lemma, we obtain that the canonical characteristic set $\D$ of $I$
is equal to the weak d-triangular subset of $\B$ of the least rank
w.r.t. $\le'$. This concludes the computation of the canonical 
characteristic set of $I$ w.r.t. the target ranking, which 
we summarize in Algorithm~\ref{a:conversionalgorithmprime}.

\begin{figure}
\begin{alg}{\sf Convert\_Prime}\label{a:conversionalgorithmprime}$\;(\C,$ $\le,$ $\le')$\\
\begin{tabular}{l}
{\sc Input:} a prime differential ideal $P = [\C]:H_\C^\infty \subset \K\{y_1,\ldots,y_n\}$\\
\hphantom{\sc Input:} with a characteristic set $\C$ w.r.t. the input ranking $\le$ \\
\hphantom{\sc Input:} with leading variables $y_1,\ldots,y_k$ and \\
\hphantom{\sc Input:} a target ranking $\le'$.\\
{\sc Output:} canonical characteristic set of $P$ w.r.t. $\le'$. \\
\end{tabular}\\
\begin{tabular}{l}
\b $M_\C := \min\left(|\C|\cdot\max\limits_{C\in\C}\ord C,\;\frac{(n-1)!}{(n-|\C|-1)!}\cdot M(\C)\right)$\\
\b $m_i := M_\C$, $1\Le i\Le k$\\
\b $\A :=$ {\sf Differentiate\&Autoreduce}$\;\left(\C,\{m_i\}_{i=1}^k\right)$\\
\b $\D :=$ {\sf Canonical\_Algebraic\_CharSet}$\;((\A):H_\A^\infty,$ $\le')$\\
\b {\bf return} {\sf minimal d-triangular subset}$\;(\D,$ $\le')$\\
\end{tabular}
\end{alg}
\end{figure}

\section{Transformation of characteristic decompositions of 
  radical differential ideals}\label{sec:characterizable}

We generalize the algebraic method for transforming characteristic 
sets of a prime differential ideal from one ranking to another to 
the case of a characterizable differential ideal. Since an ideal
characterizable w.r.t. one ranking may not be characterizable
w.r.t. another, we need to reformulate the problem: given
a characterizable differential ideal $I$ with a characteristic set
$\C$ w.r.t. a ranking $\le$, compute a characteristic decomposition
of $I$ w.r.t. another ranking $\le'$ {\it algebraically}. By analogy
with the prime case, an algebraic computation here means finding a 
sufficient differential prolongation of $\C$, which defines a
characterizable algebraic sub-ideal $\bar I$ in $I$, such that
a differential characteristic decomposition of $I$ w.r.t. $\le'$ 
can be extracted from an algebraic characteristic decomposition
of $\bar I$ w.r.t. $\le'$. 

We note that, given a characteristic
decomposition of a radical differential ideal w.r.t. one ranking,
we can obtain its characteristic decomposition w.r.t. another
ranking algebraically by solving the above problem for each
characterizable component.

All results of this section hold in the partial differential case,
except for the bound in Section~\ref{sec:primecomps}, which so
far is known only for the ordinary case. 

\subsection{Differential prolongation}

\begin{defn}
   Let $F$ be a subset in a ring $\K\{Y\}$ 
   of partial differential polynomials with a set of derivations
   $\Delta$. A set $G\subset\Theta F$ is called
   a {\it differential prolongation} of $F$, if $F\subset G$ and 
   the complement of $G$, 
   $\Theta F\setminus G$, is invariant w.r.t. differentiation, i.e.,
   for all $f\in\Theta F\setminus G$ and $\delta\in\Delta$,
   $\delta f\in\Theta F\setminus G$. 
\end{defn}   

A particular case of a differential prolongation of a weak d-triangular
set $F$ is 
$F$ itself. If $F=\C$ is autoreduced and coherent then, according to 
\cite[Lemma 6, page 137]{Kol} and \cite[Lemma 6.1 and Theorem 6.2]{Fac}, the differential ideal $I=[\C]:H_\C^\infty$
is prime, respectively characterizable iff the algebraic ideal 
$J=(\C):H_\C^\infty$ is prime, respectively characterizable. 
The ideal $J$ can be  considered either as an algebraic ideal in the ring
of differential polynomials $\K\{Y\}$ or as an ideal in the polynomial 
subring $\K[Z_\C]$, where $Z_\C=L\cup N$, $L=\ld\C$, $N=\Theta Y\setminus\Theta L$, since
the fact that $\C$ is autoreduced implies $\C\subset\K[Z_\C]$.
The Rosenfeld Lemma states that 
$$[\C]:H_\C^\infty\cap\K[Z_\C]=(\C):H_\C^\infty,$$
where the latter ideal is considered in $\K[Z_\C]$.
Moreover, a set $\D$ is a differential characteristic 
set of $I$ iff $\D$ is an algebraic characteristic set of $J$
(if the latter is considered in $\K[Z_\C]$, otherwise we
need to impose an additional requirement that $\D$ is 
differentially autoreduced).  
In particular, the canonical characteristic sets of $I$ and $J$ 
(differential and algebraic, respectively) coincide (for this
statement, it does not matter in which ring to consider $J$, since
the canonical characteristic set of an ideal is the same regardless
of the ring in which the ideal is considered). 

Now, if we consider a differential prolongation $\D$ of $\C$ and
the corresponding polynomial subring $\K[Z_\D]$, where
$Z_\D=\bar L\cup N$, $\bar L=\ld\D$, 
$N=\Theta Y\setminus\Theta L=\Theta Y\setminus\Theta\bar L$,
then $\D$ is not necessarily a subset of $\K[Z_\D]$:
\begin{exmp}\label{example:contr}
Let $\C = y',$ $x+y$ with the elimination ranking $y < x$ and a 
prolongation $$\D = y',\ x+y,\ x'+y',\ x''+y''.$$ Then 
$$\bar L = y',\ x,\ x',\ x'',\quad N = y.$$ Hence, we have that $x''+y'' \notin \K[Z_\D].$ Also, $$[\C]:H_\C^\infty\cap \K[Z_\D] = (y',\ x+y,\ x',\ x'')$$ 
and $x'' \notin (\D):H_\D^\infty.$ 
\end{exmp}
Therefore, we
need to distinguish between two ideals $I_\D:=(\D):H_\D^\infty$
in $\K\{Y\}$ and $\bar I_\D:=I\cap\K[Z_\D]$ in $\K[Z_\D]$. 




The algebraic ideal $\bar I_\D$ depends only on the set of leaders
$\bar L$ of the differential prolongation of $\C$. In other words, 
for any characterizing set $\tilde\C$ of $I$ and its differential
prolongation $\tilde\D$ with $\ld\tilde\D=\ld\D=\bar L$, we have 
$\bar I_{\tilde\D}=\bar I_\D$. We call $\bar I_{\bar L}:=\bar I_\D$
a {\it prolongation ideal} of the ideal $I$. 

Next, we study the properties of the prolongation ideals. 
The following lemma gives a criterion for a prolongation ideal
to be prime or characterizable.

\begin{lem} \label{l:prolchar}
  Let $\C$ be a coherent autoreduced set, and let $\D$ be a
  differential prolongation of $\C$. Then the differential ideal
  $1 \notin I=[\C]:H_\C^\infty$ is prime, respectively characterizable, iff 
  the corresponding prolongation ideal $\bar I_\D$ is prime, 
  respectively characterizable.
\end{lem}
\begin{pf} If $I$ is prime then its restriction $I\cap \K[Z_\D] = \bar
  I_\D$ is also prime. If $\bar I_\D$ is prime than its restriction
$\bar I_\D \cap \K[Z_\C] = (\C):H_\C^\infty$ is prime and, thus, $I$ is
prime.

Let $I$ be a characterizable differential ideal. 
We will show that set $\A$ given by
formula~\eqref{CharSetCostruction} characterizes
the prolongation ideal $\bar I_\D.$ We have 
$\bar I_\D \subset (\A):H_\A^\infty$. 
Indeed, by Lemma~\ref{InitialBelongs}, sets $\A$ and $\D$
have the same ranks, whence they have the same sets of
reduced polynomials. In particular, 
since $\D$ is a differential prolongation of the 
characteristic set $\C$, the ideal $\bar I_\D$ has no
non-zero polynomials reduced w.r.t. $\D$, and hence w.r.t. $\A$. 

Now note that $(\A):H_\A^\infty \subset I$ and
$\A \subset \K[Z_\D].$ Hence, $\bar I_\D = (\A):H_\A^\infty$ and 
$\A$ is a characteristic set of $\bar I_\D.$ Thus, $\bar I_\D$ 
is characterizable.

Since $\C\subset\K[Z_\D]$ and $(\C):H_\C^\infty=I\cap\K[Z_\C]$, we
have
$$(\C):H_\C^\infty  = (\A):H_\A^\infty\cap\K[Z_\C].$$
\end{pf}

The next lemma establishes a relation between the characteristic
sets of a characterizable differential ideal $I$ and the algebraic
characteristic sets of its prolongation ideals.
\begin{lem} \label{l:prolcharset}
  Let $\C$ be a characteristic set of the differential
  ideal $1 \notin I=[\C]:H_\C^\infty$, let $\bar L$ be a differential
  prolongation of $L=\ld\C$, and let $\bar I_{\bar L}$ be
  the corresponding prolongation ideal.

  Then a characterizing set $\A$ of $\bar I_{\bar L}$ can be obtained
  from $\C$ as
  \begin{equation}\label{CharSetCostruction}
\A:=\{\algrem(f,\B\setminus\{f\})\;|\;f\in\B,\;\ld f\in\bar L\},
\end{equation}
  where $\B$ is any triangular subset of $\Theta\C$ satisfying
  $\ld\B=\ld\Theta\C$. 

  Vice versa, given a characterizing set $\A$ of $\bar I_{\bar L}$, 
  let $\T$ be a weak d-triangular subset of $\A$ of the least
  rank. If $\T$ is differentially autoreduced, then it is a
  characterizing set of $I$. In particular, if $\A$ is the canonical
  characteristic set of $\bar I_{\bar L}$, then $\T$ is the 
  canonical characteristic set of $I$.
\end{lem}
\begin{pf}
Since $I$ is characterizable,  $\bar I_{\bar L}$ is also 
characterizable by Lemma~\ref{l:prolchar} and $\A$ is its characteristic set.
The other way follows from Lemma~\ref{l:rkweakdtr}.
\end{pf}

   In the ordinary case, the triangular set $\B$ considered in 
   the above lemma is unique. Moreover, set $\A$ can be equivalently 
   obtained as
   $$\A:={\sf Differentiate\&Autoreduce}(\C,\{m_i\}),$$ 
   where the numbers $\{m_i\}$ are the maximal orders of 
   derivatives of the leading differential indeterminates of $\C$ 
   occurring in the prolongation $\bar L$. It is preferable to 
   compute $\A$ in this way, because {\sf Differentiate\&Autoreduce}
   provides a bound on the orders of non-leading derivatives 
   occurring in $\A$, which can be used for establishing
   complexity estimates for the entire transformation algorithm.
 
   A generalization of Algorithm {\sf Differentiate\&Autoreduce} to 
   the partial case is an interesting open problem. 
   Moreover, in the partial case, there
   may be uncountably infinitely many triangular subsets of $\Theta\C$
   whose leaders coincide with $\ld\Theta\C$. Thus, not every such set can be 
   enumerated by an algorithmic procedure. However, it is easy to write
   a procedure that would enumerate a particular subset of $\Theta\C$,
   given $\C$;
   this procedure makes the computation of the set of algebraic 
   pseudo-remainders algorithmic as well. If one would like to choose
   the subset $\B$ in a systematic way, we suggest to use the ideas
   from the theory of monomial involutive divisions (see \cite{Involutive}). 

According to \cite[Theorem 4.13]{Dif}, there is a one-to-one 
correspondence between the essential prime components of 
a characterizable differential ideal $[\C]:H_\C^\infty$ and the minimal prime
components of the corresponding algebraic ideal $(\C):H_\C^\infty$. 
The following lemma generalizes this result to prolongation ideals.

\begin{lem} \label{l:proldecomp}
  Let $\C$ be a characteristic set of the differential
  ideal $I=[\C]:H_\C^\infty$, let $\bar L$ be a differential
  prolongation of $L=\ld\C$, and let $\bar I_{\bar L}$ be
  the corresponding prolongation ideal. 

  Let $I=P_1\cap\ldots\cap P_k$ be the essential prime decomposition
  of $I$, and let $(\bar P_i)_{\bar L}$ be the prolongation ideals
  corresponding to $P_i$, $i=1,\ldots,k$. Then
  $$\bar I_{\bar L}=(\bar P_1)_{\bar L}\cap\ldots\cap
  (\bar P_k)_{\bar L}$$
  is the minimal prime decomposition of $\bar I_{\bar L}$.
\end{lem}
\begin{pf} Since $\bar I_{\bar L} = I\cap \K[Z_{\bar L}]$,
$$\bar I_{\bar L} = (P_1\cap \K[Z_{\bar L}])\cap\ldots\cap (P_k\cap \K[Z_{\bar L}]) = (\bar P_1)_{\bar L}\cap\ldots\cap
  (\bar P_k)_{\bar L}$$ 
is a prime decomposition of the ideal $\bar I_{\bar L}$.
Suppose that it is not minimal. Since $(\C):H_\C^\infty=\bar I_{\bar
  L}\cap \K[Z_\C]$, 
$$(\C):H_\C^\infty=\bigl((\bar P_1)_{\bar L}\cap\K[Z_\C]\bigr)\cap\ldots\cap
                   \bigl((\bar P_k)_{\bar L}\cap\K[Z_\C]\bigr)$$
is a prime decomposition of the ideal $(\C):H_\C^\infty$, which 
is also not minimal. But the latter contradicts the fact that
$(\bar P_i)_{\bar L}\cap\K[Z_\C]=P_i\cap\K[Z_\C]$, $1\Le i\Le k$, 
and 
$$(\C):H_\C^\infty=(P_1\cap \K[Z_\C])\cap\ldots\cap
(P_k\cap\K[Z_\C])$$
is the minimal prime decomposition. 
\end{pf}

\subsection{A bound for characteristic sets of prime components} \label{sec:primecomps}

Let $I=[\C]:H_\C^\infty$ be a characterizable differential ideal with a
characteristic set $\C$ w.r.t. a ranking $\le$. Let $L=\ld_\le\C$,
and let $\bar L$ be a differential prolongation of $L$.
From the previous section we know that the prolongation ideal $\bar
I_{\bar L}$ is characterizable (Lemma~\ref{l:prolchar}) and its
minimal prime components correspond to the essential prime
components of $I$ (Lemma~\ref{l:proldecomp}). We would like to find a sufficient differential
prolongation $\bar L$ such that the minimal prime components
of $\bar I_{\bar L}$ contain differential characteristic sets
of the corresponding essential prime components of $I$ w.r.t. any
other ranking $\le'$. 

First of all, according to \cite[Theorem 4.13]{Dif}, 
a differential characteristic set of an essential prime component 
of $I$ coincides with an algebraic characteristic set of the
corresponding minimal prime component of the ideal $(\C):H_\C^\infty$.
This implies that every essential prime component $P$ of $I$ has
a characteristic set $\C_P$ satisfying the bound $m_y(\C_P)\Le m_y(\C)$
on the orders of derivatives of any differential indeterminate $y\in
Y$ occurring in $\C_P$. 

For the ordinary case, as was shown in Section~\ref{sec:primebound}, we thus have
a bound $M_\C$ on the orders of derivatives occurring in the canonical characteristic 
sets of the essential prime components of $I$ w.r.t. any other ranking
$\le'$. For the partial differential case, such a bound is not
known, but let us assume that we can compute such a bound $M_\C$
also for the partial case.\footnote{Of course, $M_\C$ can be obtained
by computing characteristic sets of the prime components w.r.t. 
the target ranking, but this would clearly defeat our purpose:  
we need a bound that can be computed from $\C$ relatively easily.} 
We need to assume that $M_\C\Ge m_y(\C)$ for all $y\in Y$. 

Let 
\begin{equation} \label{suffprol}
\bar L=\{\theta u\;|\;u\in L,\;\ord\theta u\Le M_\C\}
\end{equation}
be the differential prolongation of $L$ up to the order $M_\C$. 
According to Lemma~\ref{l:proldecomp}, the minimal prime components
of $\bar I_{\bar L}$ contain all polynomials of the corresponding
essential prime components of $I$ of order less than or equal to
$M_\C$. Thus, they also contain the canonical characteristic sets of 
the corresponding essential prime components of $I$ w.r.t. any
other ranking $\le'$. In what follows, we will denote the above 
differential prolongation $\bar I_{\bar L}$ simply by $\bar I$. 
Applying Lemma~\ref{l:prolcharset}, we compute a characteristic 
set of $\bar I$ w.r.t. $\le$. 

\subsection{Algebraic bi-characteristic decomposition}

So, we have the differential ideal $I$ which is characterizable w.r.t.
the ranking $\le$ and would like to give a characteristic decomposition
of $I$ w.r.t. $\le'.$ We have constructed the prolongation 
algebraic ideal $\bar I$ which is characterizable w.r.t. $\le$ with
a characteristic set $\A$ given by formula~\eqref{CharSetCostruction}. Let 
\begin{equation} \label{bidec}
\bar I=\bar J_1\cap\ldots\cap \bar J_k 
\end{equation}
be a bi-characteristic
decomposition of $\bar I$ w.r.t. $\le$ and $\le'$. That is, 
each component $\bar J_i$, $1\Le i\Le k$, is an algebraic ideal characterizable
w.r.t. both rankings with the canonical characteristic sets 
$\A_i$ and $\B_i$ w.r.t. $\le$ and $\le'$, respectively.

\begin{figure}
\begin{alg}{\sf Algebraic-Bicharacteristic-Decomposition}\label{BiChar}
$(\C,\le,\le')$\\
\begin{tabular}{l}
{\sc Input:} characterizing set $\C$ of a characterizable
             algebraic ideal $I$ \\
\hphantom{\sc Input:} w.r.t. an ordering $\le$ on variables \\
\hphantom{\sc Input:} and another ordering $\le'$\\
{\sc Output:} a finite set $T=\{(\C_i,\D_i)\;|\;i\in\mathfrak{I}\}$, where \\
\hphantom{\sc Output:} for every $i\in\mathfrak{I}$, $\C_i$ and $\D_i$ are algebraic 
                      characterizing sets\\
\hphantom{\sc Output:} of the same ideal $I_i$ w.r.t. $\le$ and $\le'$,
             respectively, and\\
\hphantom{\sc Output:} $I=\cap_{i\in\mathfrak{I}} I_i$\\
\\
\b $\le_s:=\le$, $\le_t:=\le'$\\
\b  ${\mathfrak C}:=\{\C\}$, $T:=\varnothing$\\
\b {\bf while} ${\mathfrak C}\neq\varnothing$ {\bf do}\\
\b \b $U:=\mathfrak{C}$, $\mathfrak{C}:=\varnothing$\\
\b \b {\bf for} $\C\in U$ {\bf do}\\
\b \b \b $J:=(\C):H_\C^\infty$ w.r.t. $\le_s$\\
\b \b \b ${\mathfrak D}:=$Algebraic-characteristic-decomposition$(\C,\le_s,\le_t)$\\
\b \b \b {\bf if} $\exists\;\D\in{\mathfrak D}$ such that
             $J=(\D):H_\D^\infty$ w.r.t. $\le_t$ {\bf then}\\
\b \b \b \b {\bf if} $\le_s=\le$ {\bf then} $T:=T\cup\{(\C,\D)\}$ {\bf
             else} $T:=T\cup\{(\D,\C)\}$\\
\b \b \b {\bf else} $\mathfrak{C}:=\mathfrak{C}\cup{\mathfrak D}$\\
\b \b \b {\bf end if}\\
\b \b {\bf end for}\\
\b \b {\bf if} $\le_s=\le$ {\bf then} $\le_s:=\le'$, $\le_t:=\le$ {\bf else}
             $\le_s:=\le$, $\le_t:=\le'$\\
\b {\bf end while}\\
\b {\bf return} $T$\\
\end{tabular}
\end{alg}
\end{figure}

Let us discuss how one can construct such a decomposition. Algorithm~\ref{BiChar} does the following. Given a characterizable
algebraic ideal $I$ with the characterizing set $\C$ w.r.t. $\le_s$, 
it first computes its (possibly redundant) algebraic 
characteristic decomposition w.r.t. 
$\le_t$ via the procedure
$$\text{\it Algebraic-characteristic-decomposition}(\C,\le_s,\le_t).$$
This procedure can be performed, for example, 
by applying the {\em Triade} algorithm \citep{MMM99}, which is
implemented in the {\tt RegularChains} library in Maple 
(see \cite{RegChains}).
A parallel implementation of this algorithm,
on a shared memory machine in {\tt Aldor} is also in progress
(see \cite{AlgDecParallel}).

If one of the characterizable components turns out to be 
equal to $I$ (note that equality of characterizable algebraic ideals
can be checked, e.g., by computing their Gr\"obner bases), then
$I$ is bi-characterizable; in this case the algorithm terminates and
outputs $T$ consisting of a single pair $(\C,\D)$ of characterizing sets of $I$ w.r.t.
$\le$ and $\le'$, respectively. If all characterizable components of $I$
contain it strictly, then, for each characterizable component, we compute its
characteristic decomposition w.r.t. $\le$ and repeat the above
strategy. 

Correctness of the algorithm follows from the 
fact that, at each iteration of the {\bf while}-loop,
$\mathfrak{C}\cup T$ provides a characteristic decomposition of $I$
w.r.t. $\le_s$ and
$T$ satisfies the requirements of the output. Termination follows
from the N\"otherian property of the polynomial ring, i.e., that every
sequence of strictly nested polynomial ideals is finite.

We note that components $\bar J_i$, for which $\ld_\le\A_i\neq\ld_\le\A$,
are redundant, i.e., they can be excluded from the right-hand side of
\eqref{bidec} without affecting the intersection. 
Indeed, if $\bar I=\bar P_1\cap\ldots\cap \bar P_l$ is the
minimal prime decomposition of $\bar I$, and 
$\bar J_i=\bar Q_{i,1}\cap\ldots\cap\bar Q_{i,l_i}$ are the minimal prime
decompositions of $\bar J_i$, $1\Le i\Le k$, then
a component $\bar J_i$ is redundant, if none of $\bar P_j$, $1\Le j\Le l$, can be found 
among $\bar Q_{i,t}$, $1\Le t\Le l_i$. But this is the case if
$\ld_\le\A_i\neq\ld_\le\A$, since by
\cite[Theorem 4.13]{Dif}
the characteristic sets of $\bar P_j$ have leaders $\ld_\le\A$, while
the characteristic sets of $\bar Q_{i,t}$ have leaders $\ld_\le\A_i$.
Therefore, we can assume that for all $1\Le i\Le k$,
$\ld_\le\A_i=\ld_\le\A$.

We prove then that {\it every} minimal prime component of $\bar J_i$ is a minimal 
prime component of $\bar I$. Indeed, every $\bar Q_{i,t}$ is a prime ideal
containing $\bar I$. Suppose that $\bar Q_{i,t}$ is not minimal, i.e.,
there is a minimal prime component $\bar P_j$ of $\bar I$ such that 
$\bar P_j\subsetneq \bar Q_{i,t}$. But the latter strict inclusion is
impossible
according to the following Lemma~\ref{l:equalideals} and Remark~\ref{r:nodifferentiations}.

\begin{lem}\label{l:equalideals}
  Let $P$ and $Q$ be two prime 
differential 
ideals whose characteristic
  sets w.r.t. $\le$ have the same sets of leaders
  Then $P\subseteq Q$ implies $P=Q$.
\end{lem}
\begin{pf} Let $\C_1$ and $\C_2$ be these characteristic sets. We have
$P = [\C_1]:H_{\C_1}^\infty$ and $Q = [\C_2]:H_{\C_2}^\infty.$ Consider
the restricted ideals $\p = (\C_1):H_{\C_1}^\infty$ and $\q = (\C_2):H_{\C_2}^\infty$ in the N\"otherian ring $\K[L,N(\C_1,\C_2)],$
where $N(\C_1,\C_2)$ is the set of non-leading variables appearing in
both $\C_1$ and $\C_2.$ From \cite[Theorem 3.2]{Fac} it follows
that both $\p$ and $\q$ are of dimension $|N(\C_1,\C_2)|.$ 

Take any $f \in \p$. It is partially reduced w.r.t. both $\C_1$ and
$\C_2$ (which are coherent and autoreduced) and belongs to $P\subset Q$. By the Rosenfeld lemma $f \in \q.$ Hence, $\p \subset \q$  and they are prime and must be equal then, because their Krull dimensions
are equal to the same number $|N(\C_1,\C_2)|.$ Hence, we have $\C_1 \subset Q$ and 
$\C_2 \subset \p \subset P$ at the same time.
Thus, according to \cite[Theorem 9]{canonical} we finally obtain that $P=Q.$
\end{pf}

\begin{rem}\label{r:nodifferentiations}
 In the above lemma, one can assume that the set of derivations is
  empty, hence the statement also holds for algebraic ideals.
\end{rem}


To summarize, for every bi-characterizable component $\bar J_i$, there
exists a subset $T_i\subset\{1,\ldots,l\}$ such that 
$$\bar J_i=\bigcap_{j\in T_i}\bar P_j$$ 
is the minimal prime decomposition of $\bar J_i$. Moreover, equality \eqref{bidec}
implies that 
$$\bigcup_{i=1}^l T_i=\{1,\ldots,l\}.$$

\subsection{Constructing differential characterizable components from 
            the algebraic ones}

Fix any of the above algebraic bi-characterizable components $\bar J=\bar
J_i$, where $1\Le i\Le k$; we have a set of indices
$T=T_i\subset\{1,\ldots,l\}$ such that 
$$\bar J=\bigcap_{j\in T}\bar P_j.$$
As above, let $\A=\A_i$ and $\B=\B_i$ be the canonical characteristic
sets of $\bar J$ w.r.t. $\le$ and $\le'$, respectively. 

According to Lemma~\ref{l:proldecomp}, each minimal prime component
$\bar P_j$ of $\bar I$ is a prolongation ideal of the corresponding essential 
prime component $P_j$ of $I$, i.e., 
$$\bar P_j=P_j\cap\K[\bar L\cup N],$$
where $I=\bigcap_{j=1}^l P_j$ is the essential prime decomposition of
$I$. 
Since $\B$ is a characterizing set of $\bar J$ w.r.t. $\le'$, the initials and separants
of $\B$ w.r.t. $\le'$ are not zero-divisors modulo $\bar J$, i.e., they do not 
belong to the minimal prime components $\bar P_j$, $j\in T$.
Since $\B$, as well as $H_{\B}$, is a subset of $\K[\bar L\cup N]$, we have therefore
$H_{\B}\cap P_j=\varnothing$, $j\in T$. 

Let $\T\subset\B$ be the weak d-triangular subset of $\B$ of the 
least rank w.r.t. $\le'$. Since $H_{\T}\subset H_{\B}$, we also have
$H_{\T}\cap P_j=\varnothing$, $j\in T$. Thus, 
$[\T]:H_{\T}^\infty\subset P_j$, $j\in T$. In particular, this
implies that $[\T]:H_{\T}^\infty\neq (1)$.

Let $\D$ be the result of differential autoreduction of $\T$ w.r.t. $\le'$,
i.e.,
$$\D=\{\drem(f,\T\setminus\{f\})\;|\;f\in \T\}.$$
Set $\D$ is differentially autoreduced. 
By definition of differential remainder, $\D\subset [\T]$. 
By Lemma~\ref{InitialBelongs}, since $[\T]:H_{\T}^\infty\neq(1)$,
we have $\rk_{\le'}\D=\rk_{\le'}\T$ and, moreover, 
$H_{\D}\subset H_{\T}^\infty+[\T]$. Therefore, 
\begin{equation} \label{incld}
[\D]:H_{\D}^\infty\subset [\T]:H_{\T}^\infty\subset P_j,\;\;j\in T.
\end{equation}

We will show that $\D$ is a characteristic set of the ideal 
$[\D]:H_{\D}^\infty$ w.r.t. $\le'$ 
by proving that every polynomial 
in the intersection $\bigcap_{j\in T_i} P_j$ reduces w.r.t. $\D$ to
zero. Given \eqref{incld}, this will also imply that 
\begin{equation} \label{inters}
[\D]:H_{\D}^\infty=\bigcap_{j\in T}P_j.
\end{equation}

Take any polynomial $f\in\bigcap_{j\in T}P_j$, and let 
$\bar f=\drem(f,\D)$, where the pseudo-remainder is computed
w.r.t. $\le'$. Since $\D\subset[\T]\subset P_j$, $j\in T$, 
we have \mbox{$\bar f\in\bigcap_{j\in T} P_j$}. 

Let $\F_j$ be the canonical characteristic set of $P_j$ w.r.t. $\le'$,
and let $\bar\F_j$ be the canonical algebraic characteristic set
of the corresponding prolongation ideal $\bar P_j$. 
We have shown in Section~\ref{sec:primecomps} that $\bar P_j$
contains $\F_j$. Thus, from Lemma~\ref{l:algdifcs} it follows that 
$\F_j$ is the weak d-triangular subset of $\bar\F_j$ of the least rank 
w.r.t. $\le'$. 
On the other hand, since $\bar P_j$ is a minimal prime component
of $\bar J$, according to \cite[Theorem 4.13]{Dif}, 
$\ld_{\le'}\bar \F_j=\ld_{\le'}\B$. This implies that 
$\ld_{\le'} \F_j=\ld_{\le'}\T=\ld_{\le'}\D$. That is,
the fact that $\bar f$ is reduced w.r.t. $\D$ implies that
it is partially reduced w.r.t. $\F_j$. 

By the Rosenfeld Lemma, 
$$\bar f\in (\F_j):H_{\F_j}^\infty\subset (\bar\F_j):H_{\bar\F_j}^\infty=\bar P_j,\;\;j\in T$$
i.e., $\bar f\in \bar J$. Now, the fact that $\bar f$ is reduced
w.r.t. $\D$ implies that it is algebraically reduced w.r.t. $\B.$
Since the latter is a characteristic set of $\bar J$, we obtain
$\bar f=0$ and the required equality \eqref{inters}.

Now we see that the ideal $[\D]:H_\D^\infty$ is characterizable
w.r.t. $\le'$. The canonical characteristic set of this ideal 
w.r.t. $\le'$ is contained in each minimal prime component
of the ideal $(\D):H_\D^\infty$, therefore it is also contained
in every $\bar P_j$, $j\in T$, and hence in $\bar J$. The ideal
$\bar J$ is contained in $[\D]:H_\D^\infty$. Thus, by
Lemma~\ref{l:algdifcs}, the canonical characteristic set of 
$[\D]:H_\D^\infty$ is equal to the weak d-triangular subset
of $\B$ of the least rank w.r.t. $\le'$.  That is, we have $$\D=\T$$
which is (w.r.t. the ranking $\le'$) the canonical characteristic set of the characterizable differential
ideal
$$[\D]:H_\D^\infty.$$ 

\subsection{The final characteristic decomposition}

In the previous section, we have shown that for each
bi-characterizable component $\bar J_i$, $1\le i\le l$, of $\bar I$ with 
the canonical characteristic set $\B_i$ w.r.t. $\le'$, 
if $\D_i$ is the weak d-triangular subset of $\B_i$ of the least rank,
then it is the canonical characteristic set of the ideal 
$[\D_i]:H_{\D_i}^\infty$. 
We have also shown that 
$$[\D_i]:H_{\D_i}^\infty=\bigcap_{j\in T_i} P_j.$$
Thus, since $\bigcup_{i=1}^l T_i=\{1,\ldots,l\}$, the following intersection
$$\bigcap_{i=1}^l[\D_i]:H_{\D_i}^\infty$$
is a characteristic decomposition of $I=P_1\cap\ldots\cap P_l$
w.r.t. $\le'$.
This concludes the algebraic computation of a characteristic
decomposition of $I$ w.r.t. the target ranking, which we summarize in
the Algorithm~\ref{a:conversionalgorithm}.

\begin{figure}
\begin{alg}{\sf Convert\_Characterizable}\label{a:conversionalgorithm}$\;(\C,$ $\le,$ $\le')$\\
\begin{tabular}{l}
{\sc Input:} set $\C$ which characterizes the ideal $[\C]:H_\C^\infty$
             w.r.t. the input ranking $\le$\\
\hphantom{\sc Input:}  and has leading variables $y_1,\ldots,y_k$ and a target ranking $\le'$.\\
{\sc Output:} characteristic decomposition of $[\C]:H_\C^\infty$ w.r.t. $\le'$. \\
\end{tabular}
\begin{tabular}{l}
\b $M_\C := \min\left(|\C|\cdot\max\limits_{C\in\C}\ord C,\;\frac{(n-1)!}{(n-|\C|-1)!}\cdot M(\C)\right)$\\
\b $m_i := M_\C$, $1\Le i\Le k$\\
\b $\A :=$ {\sf Differentiate\&Autoreduce}$\;\left(\C,\{m_i\}_{i=1}^k\right)$\\
\b $\mathfrak{D} :=$ {\sf Bi-characterizable\_Canonical\_Decomposition}\;$((\A):H_\A^\infty,$ $\le,$ $\le')$\\
\b $\mathfrak{C} := \{{\sf minimal\ d-triangular\ subset}\;(\D, \le')\;|\;\D\in\mathfrak{D}\}$\\
\b {\bf return} $\mathfrak{C}$
\end{tabular}
\end{alg}
\end{figure}

Now, in order to convert a characteristic decomposition 
$$I = \bigcap_{i=1}^p[\C_i]:H_{\C_i}^\infty$$ of a radical differential
ideal $I$ w.r.t. $\le$ to a ranking $\le'$, one just applies 
Algorithm~\ref{a:conversionalgorithm} to each characterizable
component $[\C_i]:H_{\C_i}^\infty$ and then collects all the results 
together in a single intersection.

\section{Conclusions}
By estimating the orders of derivatives, we have shown that,
given a set of ordinary differential polynomials specifying a radical
differential ideal $I$, one can construct a N\"otherian ring of algebraic
polynomials, in which the computation of a characteristic
decomposition of $I$ is actually performed. This does not
mean that the computation is completely algebraic: differentiations
are allowed, but they never lead out of the constructed algebraic ring.

For the problem of converting a characteristic decomposition
of a radical differential ideal from one ranking to another, we
have proposed an algorithm, which first differentiates the input
polynomials sufficiently many times, and then performs the
conversion completely algebraically, without using differentiation
at all. The algorithm is applicable in the partial differential case,
but the bound for the number of differentiations of the input
polynomials is given for the ordinary case only.

We conjecture that, if one can solve the first problem
of computing a characteristic decomposition of a radical
differential ideal from generators completely algebraically,
i.e., by an algorithm that first differentiates the input
polynomials sufficiently many times, and then computes
the decomposition without using differentiations, then
one can also solve the Ritt problem of computing
an irredundant prime (or characteristic) decomposition
of a radical differential ideal.

\begin{ack} We thank Michael F. Singer, Fran\c{c}ois Boulier,
William Sit, \'Evelyne Hubert, Evgeniy Pankratiev, and the 
referees for their important suggestions.
\end{ack}

\bibliographystyle{elsart-harv}
\bibliography{rg}

\end{document}